\documentclass[mnsc,nonblindrev]{informs3}

\OneAndAHalfSpacedXI

\usepackage{natbib}
 \bibpunct[, ]{(}{)}{,}{a}{}{,}
 \def\bibfont{\small}

\TheoremsNumberedThrough     

\EquationsNumberedThrough

\MANUSCRIPTNO{}

\usepackage{amsmath, amssymb, bbold, mathtools}
\usepackage{footnote}
\usepackage{acronym}
\usepackage{subfigure}
\usepackage{import}
\usepackage{cite}
\usepackage{enumerate}

\acrodef{dro}[DRO]{distributionally robust optimization}
\acrodef{ldt}[LDT]{large deviation theory}
\acrodef{ldp}[LDP]{large deviation principle}
\acrodef{lln}[LLN]{law of large numbers}
\acrodef{kl}[KL]{Kullback-Leibler}
\acrodef{iid}[{i.i.d.\ \!\!}]{independent identically distributed}
\acrodef{qp}[QP]{quadratic program}
\acrodef{qcqp}[QCQP]{quadratically constrained quadratic program}
\acrodef{vod}[VoD]{value of data}
\acrodef{saa}[SAA]{stochastic average approximation}
\acrodef{fl}[FL]{Fenchel-Legendre}

\usepackage{tikz}
\usepackage{pgfplots}
\usetikzlibrary{plotmarks,arrows,shapes,backgrounds}
\usetikzlibrary{automata}
\usetikzlibrary{matrix,positioning}
\usetikzlibrary{decorations.pathmorphing,calc}
\usepackage{tikz-3dplot}
\usetikzlibrary{external}
\newcommand{\norm}[1]{\left\|#1\right\|}
\newcommand{\abs}[1]{\left|#1\right|}

\newcommand{\D}[2]{\mathrm I(#1 , #2 )}

\newcommand{\set}[2]{\left\{ #1\ : \ #2 \right\}}

\renewcommand{\mc}{\mathcal}
\newcommand{\mb}{\mathbb}

\def\d{\mathrm{d}}
\def\st{\mathrm{s.t.}}

\DeclareMathOperator{\cl}{cl}
\DeclareMathOperator{\interior}{int}

\begin{document}

\RUNAUTHOR{Van Parys et al.}

\RUNTITLE{Distributionally Robust Optimization is Optimal}

\TITLE{From Data to Decisions: \\ Distributionally Robust Optimization is Optimal}

\ARTICLEAUTHORS{%
\AUTHOR{Bart P.G.\ Van Parys}
\AFF{Operations Research Center, Massachusetts Institute of Technology, \EMAIL{vanparys@mit.edu}}
\AUTHOR{Peyman Mohajerin Esfahani}
\AFF{Delft Center for Systems and Control, Technische Universiteit Delft, \EMAIL{p.mohajerinesfahani@tudelft.nl}}
\AUTHOR{Daniel Kuhn}
\AFF{Risk Analytics and Optimization Chair, Ecole Polytechnique F\'ed\'erale de Lausanne, \EMAIL{daniel.kuhn@epfl.ch}}

}

\ABSTRACT{%
We study stochastic programs where the decision-maker cannot observe the distribution of the exogenous uncertainties but has access to a finite set of independent samples from this distribution. In this setting, the goal is to find a procedure that transforms the data to an estimate of the expected cost function under the unknown data-generating distribution, {\em i.e.}, a {\em predictor}, and an optimizer of the estimated cost function that serves as a near-optimal candidate decision, {\em i.e.}, a {\em prescriptor}. As functions of the data, predictors and prescriptors constitute statistical estimators. We propose a meta-optimization problem to find the least conservative predictors and prescriptors subject to constraints on their {\em out-of-sample disappointment}. The out-of-sample disappointment quantifies the probability that the actual expected cost of the candidate decision under the unknown true distribution exceeds its predicted cost. Leveraging tools from large deviations theory, we prove that this meta-optimization problem admits a unique solution: The best predictor-prescriptor-pair is obtained by solving a distributionally robust optimization problem over all distributions within a given relative entropy distance from the empirical distribution of the data.
}%

\KEYWORDS{Data-Driven Optimization, Distributionally Robust Optimization, Large Deviations Theory, Relative Entropy, Convex Optimization, Observed Fisher Information} \HISTORY{}

\maketitle

\acresetall

\section{Introduction}
We study static decision problems under uncertainty, where the decision maker cannot observe the probability distribution of the uncertain problem parameters but has access to a finite number of independent samples from this distribution. Classical stochastic programming uses this data only indirectly. The data serves as the input for a statistical estimation problem that aims to infer the distribution of the uncertain problem parameters. The estimated distribution then serves as an input for an optimization problem that outputs a near-optimal decision as well as an estimate of the expected cost incurred by this decision. Thus, classical stochastic programming separates the decision-making process into an estimation phase and a subsequent optimization phase. The estimation method is typically selected with the goal to achieve maximum prediction accuracy but without tailoring it to the optimization problem at hand.

In this paper we develop a method of data-driven stochastic programming that avoids the artificial decoupling of estimation and optimization and that chooses an estimator that adapts to the underlying optimization problem. Specifically, we model data-driven solutions to a stochastic program through a {\em predictor} and its corresponding {\em prescriptor}. For any fixed feasible decision, the predictor maps the observable data to an estimate of the decision's expected cost. The prescriptor, on the other hand, computes a decision that minimizes the cost estimated by the predictor.

The set of all possible predictors and their induced prescriptors is vast. Indeed, there are countless possibilities to estimate the expected costs of a fixed decision from data, {\em e.g.}, via the popular sample average approximation \citep[Chapter~5]{shapiro2014lectures}, by postulating a parametric model for the exogenous uncertainties and estimating its parameters via maximum likelihood estimation \citep{dupacovawets1988}, or through kernel density estimation \citep{parpas2015}. Recently, it has become fashionable to construct conservative (pessimistic) estimates of the expected costs via methods of distributionally robust optimization. In this setting, the available data is used to generate an {\em ambiguity set} that represents a confidence region in the space of probability distributions and contains the unknown data-generating distribution with high probability. The expected cost of a fixed decision under the unknown true distribution is then estimated by the worst-case expectation over all distributions in the ambiguity set. Since the ambiguity set constitutes a confidence region for the unknown true distribution, the worst-case expectation represents an upper confidence bound on the true expected cost. The ambiguity set can be defined, for example, through confidence intervals for the distribution's moments \citep{delage2010distributionally}. Alternatively, the ambiguity set may contain all distributions that achieve a prescribed level of likelihood \citep{wang2009likelihood}, that pass a statistical hypothesis test \citep{bertsimas2014robust} or that are sufficiently close to a reference distribution with respect to a probability metric such as the Prokhorov metric \citep{erdogan2006}, the Wasserstein distance \citep{pflug2007,esfahani2015data,ZG15:Wasserstein}, the total variation distance \citep{sun2016} or the $L^1$-norm \citep{jiang2015}. \citet{ben2013robust} have shown that confidence sets for distributions can also be constructed using $\phi$-divergences such as the Pearson divergence, the Burg entropy or the Kullback-Leibler divergence. More recently, \citet{bayraksanlove2015} provide a systematic classification of $\phi$-divergences and investigate the richness of the corresponding ambiguity sets. %\citep{shapiro2015}

Given the numerous possibilities for constructing predictors from a given dataset, it is easy to lose oversight. In practice, predictors are often selected manually from within a small menu with the goal to meet certain statistical and/or computational requirements. However, there are typically many different predictors that exhibit the desired properties, and there always remains some doubt as to whether the chosen predictor is best suited for the particular decision problem at hand. In this paper we propose a principled approach to data-driven stochastic programming by solving a meta-optimization problem over a rich class of predictor-prescriptor-pairs including, among others, all examples reviewed above. This meta-optimization problem aims to find the least conservative ({\em i.e.}, pointwise smallest) prescriptor whose {\em out-of-sample disappointment} decays at a prescribed exponential rate~$r$ as the sample size tends to infinity---irrespective of the true data-generating distribution. The out-of-sample disappointment quantifies the probability that the {\em actual} expected cost of the prescriptor exceeds its {\em predicted} cost. Put differently, it represents the probability that the predicted cost of a candidate decision is over-optimistic and leads to disappointment in out-of-sample tests. Thus, the proposed meta-optimization problem tries to identify the predictor-prescriptor-pairs that overestimate the expected out-of-sample costs by the least amount possible without risking disappointment under {\em any} thinkable data-generating distribution.

Our main results can be summarized as follows.
\begin{itemize}
\item By leveraging Sanov's theorem from large deviations theory, we prove that the meta-optimization problem admits a unique optimal solution for any given stochastic program. 
\item We show that the optimal data-driven predictor estimates the expected costs under the unknown true distribution by a worst-case expectation over all distributions within a given relative entropy distance from the empirical distribution of the data. This suggests that, among all possible data-driven solutions, a distributionally robust approach based on a relative entropy ambiguity set is optimal. This is perhaps surprising because the meta-optimization problem does not impose any structure on the predictors, which are generic functions of the data. In particular, there is no requirement forcing predictors to admit a distributionally robust interpretation. 
\item {In contrast to most of the existing work on data-driven distributionally robust optimization, our relative entropy ambiguity set does {\em not} play the role of a confidence region that contains the unknown data-generating distribution with a prescribed level of probability (see the discussions of \citep{lam2016robust,gupta2015near} below for exceptions).} Instead, the radius of the relative entropy ambiguity set coincides with the desired exponential decay rate~$r$ of the out-of-sample disappointment imposed by the meta-optimization problem.
\item We prove that the optimal (distributionally robust) predictor admits a dual representation as the optimal value of a one-dimensional convex optimization problem that can be solved highly efficiently. For continuously distributed problem parameters this representation seems to be new.
\end{itemize}

To our best knowledge, we are the first to recognize the optimality of distributionally robust optimization in its ability to transform data to predictors and prescriptors. The optimal distributionally robust predictor identified in this paper can be evaluated by solving a tractable convex optimization problem. Under standard convexity assumptions about the feasible set and the cost function of the stochastic program, the corresponding optimal prescriptor can also be evaluated in polynomial time. Although perhaps desirable, the tractability and distributionally robust nature of the optimal predictor-prescriptor-pair are not dictated {\em ex ante} but emerge naturally. 

Relative entropy ambiguity sets have already attracted considerable interest in distributionally robust optimization \citep{ben2013robust, calafiore2007ambiguous, hu2013kullback, lam2016robust, wang2009likelihood}. Note, however, that the relative entropy constitutes an {\em asymmetric} distance measure between two distributions. The asymmetry implies, among others, that the first distribution must be absolutely continuous to the second one but not {\em vice versa}. Thus, ambiguity sets can be constructed in two different ways by designating the reference distribution either as the first or as the {\em second} argument of the relative entropy.  All papers listed above favor the second option, and thus the emerging ambiguity sets contain only distributions that are absolutely continuous to the reference distribution. Maybe surprisingly, the optimal predictor resulting from our meta-optimization problem uses the reference distribution as the {\em first} argument of the relative entropy instead. Thus, the reference distribution is absolutely continuous to every distribution in the emerging ambiguity set. Relative entropy balls of this kind have previously been studied by \citet{gupta2015near}, \citet{lam2016} { and \citet{bertsimas2018data}}.

Adopting a Bayesian perspective, \citet{gupta2015near} determines the smallest ambiguity sets that contain the unknown data-generating distribution with a prescribed level of confidence as the sample size tends to infinity. Both Pearson divergence and relative entropy ambiguity sets with properly scaled radii are optimal in this setting. In the terminology of the present paper, \citet{gupta2015near} thus restricts attention to the subclass of distributionally robust predictors and operates with an asymptotic notion of optimality. The meta-optimization problem proposed here entails a stronger notion of optimality, under which the distributionally robust predictor with relative entropy ambiguity set emerges as the unique optimizer. \citet{lam2016} also seeks distributionally robust predictors that trade conservatism for out-of-sample performance. He studies the probability that the estimated expected cost function dominates the actual expected cost function uniformly across all decisions, and he calls a predictor optimal if this probability is asymptotically equal to a prescribed confidence level. Using the empirical likelihood theorem of \citet{owen1988empirical}, he shows that Pearson divergence and relative entropy ambiguity sets with properly scaled radii are optimal in this sense. This notion of optimality has again an asymptotic flavor {in the sense that it refers to {\em sequences} of ambiguity sets that converge to a singleton}, and it admits multiple optimizers.

The rest of the paper unfolds as follows. Section~\ref{sec:data-driven} provides a formal introduction to data-driven stochastic programming on finite state spaces and develops the meta-optimization problem for identifying the best predictor-prescriptor-pair. Section~\ref{sec:ldt} reviews weak and strong large deviation principles, which are then used in Section~\ref{sec:distr_rob_optimization} to determine the unique optimal solution of the meta-optimization problem. An extension to continuous state spaces is discussed in Section~\ref{sec:continuous}.

\paragraph{\bf Notation:} The natural logarithm of $p\in\Re_+$ is denoted by $\log(p)$, where we use the conventions $0 \log(0/p)=0$ for any $p\geq 0$ and $p'\log(p'/0) =\infty$ for any $p'>0$. 
A function $f:\mc P\to X$ from $\mc P\subseteq \Re^d$ to $X\subseteq \Re^n$ is called quasi-continuous at $\mb P\in \mc P$ if for every $\epsilon>0$ and neighborhood $U\subseteq \mc P$ of $\mb P$ there is a non-empty open set $V\subseteq U$ with $|f(\mb P)-f(\mb Q)|\leq \epsilon$ for all $\mb Q\in V$. Note that $V$ does not necessarily contain $\mb P$. For any logical statement~$\mc E$, the indicator function $\mb 1_\mc E$ evaluates to~1 if $\mc E$ is true and to $0$ otherwise.

%%%%%%%%%%%%%%%%% 
%% SEC. Stochastic optimization
%%%%%%%%%%%%%%%%% 

\section{Data-driven stochastic programming}
\label{sec:data-driven}

Stochastic programming is a powerful modeling paradigm for taking informed decisions in an uncertain environment. A generic single-stage stochastic program can be represented as
\begin{equation}
  \label{eq:prescription_problem}
  \mathop{\text{minimize}}_{x \in X} \,\mb E_{\mb P^\star} \left[ \gamma(x, \xi)\right].
\end{equation}
Here, the goal is to minimize the expected value of a cost function $\gamma(x,\xi)\in\Re$, which depends both on a decision variable $x\in X$ and a random parameter $\xi\in \Xi$ governed by a probability distribution~$\mb P^\star$. We will assume that the cost $\gamma(x,\xi)$ is continuous in $x$ for every fixed $\xi\in\Xi$, the feasible set $X\subseteq \Re^n$ is compact, and $\Xi=\{1,\ldots,d\}$ is finite. Thus, $\xi$ has $d$ distinct scenarios that are represented---without loss of generality---by the integers $1,\ldots,d$. We will relax this requirement in Section~\ref{sec:continuous}, where $\Xi$ will be modeled as an arbitrary compact subset of $\Re^d$.  A wide spectrum of decision problems can be cast as instances of~\eqref{eq:prescription_problem}. \citet{shapiro2014lectures} point out, for example, that~\eqref{eq:prescription_problem} can be viewed as the first stage of a two-stage stochastic program, where the cost function $\gamma(x, \xi)$ embodies the optimal value of a subordinate second-stage problem. Alternatively, problem~\eqref{eq:prescription_problem} may also be interpreted as a generic learning problem in the spirit of statistical learning theory.

In the following, we distinguish the {\em prediction problem}, which merely aims to predict the expected cost associated with a fixed decision $x$, and the {\em prescription problem}, which seeks to identify a decision $x^\star$ that minimizes the expected cost across all $x\in X$. 

Any attempt to solve the prescription problem seems futile unless there is a procedure for solving the corresponding prediction problem. The generic prediction problem is closely related to what \citet{lemaitre2010introduction} call an uncertainty quantification problem and is therefore of prime interest in its own right. Throughout the rest of the paper, we thus analyze prediction and prescription problems on equal footing. 

In the what follows we formalize the notion of a data-driven solution to the prescription and prediction problems, respectively. Furthermore, we introduce the basic assumptions as well as the notation used throughout the remainder of the paper.

\subsection{Data-driven predictors and prescriptors}
\label{ssec:dd_predictions_prescriptions}

If the distribution $\mb P^\star$ of $\xi$ is unobservable and must be estimated from a training dataset consisting of finitely many independent samples from $\mb P^\star$, we lack essential information to evaluate the expected cost of any fixed decision and to solve the stochastic program~\eqref{eq:prescription_problem}. The standard approach to overcome this deficiency is to approximate $\mb P^\star$ with a parametric or non-parametric estimate $\hat {\mb P}$ inferred from the samples and to minimize the expected cost under $\hat {\mb P}$ instead of the true expected cost under $\mb P^\star$. However, if we calibrate a stochastic program to a training data set and evaluate its optimal decision on a test data set, then the resulting test performance is often disappointing---even if the two datasets are sampled independently from $\mb P^\star$. This phenomenon has been observed in many different contexts. It is particularly pronounced in finance, where \citet{michaud1989markowitz} refers to it as the `error maximization effect' of portfolio optimization, and in statistics or machine learning, where it is known as `overfitting'. In decision analysis, \citet{smith2006optimizer} refer to it as the `optimizer's curse'. Thus, when working with data instead of exact probability distributions, one should safeguard against solutions that display promising in-sample performance but lead to out-of-sample disappointment.

Initially the distribution $\mb P^\star$ is only known to belong to the probability simplex $\mc P=\{\mb P\in\Re_+^d:\sum_{i\in\Xi} \mb P(i)=1\}$. Over time, however, independent samples $\xi_t$, $t\in\mb N$, from $\mb P^\star$ are revealed to the decision maker that provide increasingly reliable statistical information about~$\mb P^\star$. 

Any $\mb P\in \mc P$ encodes a possible probabilistic model for the data process. Thus, by slight abuse of terminology, we will henceforth refer to the distributions $\mb P\in \mc P$ as models and to $ \mc P$ as the model class. Evidently, the true model $\mb P^\star$ is an (albeit unknown) element of~$ \mc P$. Next, we introduce model-based predictors and prescriptors corresponding to the stochastic program \eqref{eq:prescription_problem}, where the true unknown distribution $\mb P^\star$ is replaced with a hypothetical model $\mb P\in \mc P$.

\begin{definition}[model-based predictors and prescriptors]
  \label{def:parametric-predictor}
  For any fixed model $\mb P\in \mc P$, we define the {\em model-based predictor} $c(x,\mb P)= \mb E_{\mb P}[\gamma(x,\xi)]= \sum_{i\in\Xi} \mb P(i)\, \gamma (x,i)$ as the expected cost of a given decision $x\in X$ and the {\em model-based prescriptor} $x^\star(\mb P)\in\arg\min_{x\in X}c(x,\mb P)$ as a decision that minimizes $c(x,\mb P)$ over $x\in X$. 
\end{definition}

Note that the model-based predictor $c(x,\mb P)$ is jointly continuous in $x$ and $\mb P$ because $\Xi$ is finite and $\gamma(x,\xi)$ is continuous in $x$ for every fixed $\xi\in\Xi$. The continuity of $c(x,\mb P)$ then guarantees via the compactness of $X$ that the model-based prescriptor $x^\star(\mb P)$ exists for every model $\mb P\in\mc P$. In view of Definition~\ref{def:parametric-predictor}, the stochastic program~\eqref{eq:prescription_problem} can be identified with the {\em prescription problem} of computing $x^\star(\mb P^\star)$. Similarly, the evaluation of the expected cost of a given decision $x\in X$ in \eqref{eq:prescription_problem} can be identified with the {\em prediction problem} of computing $c(x,\mb P^\star)$. These prediction and prescription problems cannot be solved, however, as they depend on the unknown true model $\mb P^\star$. 

If one has only access to a finite set $\{\xi_t\}_{t=1}^T$ of independent samples from $\mb P^\star$ instead of $\mb P^\star$ itself, then it may be useful to construct an empirical estimator for~$\mb P^\star$.

\begin{definition}[Empirical distribution]
\label{def:empirical}
  The empirical distribution $\hat {\mb P}_T$ corresponding to the sample path $\{\xi_t\}_{t=1}^T$ of length $T$ is defined through
    \[
      \hat {\mb P}_T(i)= \frac1T \sum_{t=1}^{T} \mb 1_{\xi_t=i} \quad \forall i\in\Xi.
    \] 
\end{definition}

Note that $\hat {\mb P}_T$ can be viewed as the vector of empirical state frequencies. Indeed, its $i^{\text{th}}$ entry records the proportion of time that the sample path spends in state~$i$.
As the samples are drawn independently, the state frequencies capture all useful statistical information about $\mb P^\star$ that can possibly be extracted from a given sample path. Note also that $\hat {\mb P}_T$ is in fact the maximum likelihood estimator of $\mb P^\star$. In the following, we will therefore approximate the unknown predictor $c(x,\mb P^\star)$ as well as the unknown prescriptor $x^\star(\mb P^\star)$ by suitable functions of the empirical distribution $\hat {\mb P}_T$.

\begin{definition}[Data-driven predictors and prescriptors]
  \label{def:dd_prediction}
  A continuous function $\hat c:X\times\mc P\to\Re$ is called a {\em data-driven predictor} if $\hat c(x,\hat {\mb P}_T)$ is used as an approximation for $c(x,\mb P^\star)$. 
  A quasi-continuous function $\hat x: \mc P\to X$ is called a {\em data-driven prescriptor} if there exists a data-driven predictor $\hat c$ with
  \[
    \hat x({\mb P'}) \in \arg\min_{x\in X}\, \hat c(x, {\mb P'})
  \]
  for all possible estimator realizations ${\mb P'}\in\mc P$, and $\hat x(\hat {\mb P}_T)$ is used as an approximation for $x^\star(\mb P^\star)$.
\end{definition}

Every data-driven predictor $\hat c$ induces a data-driven prescriptor $\hat x$. To see this, note that the `$\arg\min$' mapping is non-empty-valued and upper semicontinuous due to Berge's maximum theorem \citep[pp.~115--116]{berge1963topological}, which applies because $\hat c$ is continuous and $X$ is both compact and independent of $\mb P'$. Corollary~4 in~\citep{Matejdes:87}, which applies because $\mc P$ is a Baire space and $X$ is a metric space, thus ensures that the `$\arg\min$' mapping admits a quasi-continuous selector, which serves as a valid data-driven prescriptor. One can show that the set of points where this quasi-continuous prescriptor is discontinuous is a meagre subset of $\mc P$ \citep{Bledsoe1952}. By the Baire category theorem, the points of continuity of the data-driven prescriptor at hand are thus dense in $\mc P$ \citep{baire1899fonctions}. Thus, data-driven prescriptors in the sense of Definition~\ref{def:dd_prediction} are `mostly'~continuous.

\begin{example}[Sample average predictor]
\label{ex:naive-predictor}
The model-based predictor $c$ introduced in Definition \ref{def:parametric-predictor} constitutes a simple data-driven predictor $\hat c=c$, that is, $c(x,\hat {\mb P}_T)$ can readily be used as a na\"ive approximation for $c(x,\mb P^\star)$. Note that the model-based predictor $c$ is indeed continuous as desired. By the definition of the empirical estimator, this na\"ive predictor approximates $c(x,\mb P^\star)$~with
\[
	c(x,\hat {\mb P}_T) = \frac{1}{T}\sum_{t=1}^T \gamma(x, \xi_t),
\]
which is readily recognized as the popular sample average approximation.
\end{example}

\subsection{Optimizing over all data-driven predictors and prescriptors}

The estimates $\hat c(x,\hat {\mb P}_T)$ and $\hat x(\hat {\mb P}_T)$ inherit the randomness from the empirical estimator~$\hat{\mb P}_T$, which is constructed from the (random) samples $\{\xi_t\}_{t=1}^T$. Note that the prediction and prescription problems are naturally interpreted as instances of statistical estimation problems. Indeed, data-driven prediction aims to estimate the expected cost $c(x, \mb P^\star)$ from data. Standard statistical estimation theory would typically endeavor to find a data-driven predictor $\hat c$ that (approximately) minimizes the mean squared error
\begin{equation*}
  \mb E \left[ |c(x, \mb P^\star) - \hat c(x, \hat {\mb P}_T)|^2 \right]
\end{equation*}
over some appropriately chosen class of predictors $\hat c$, where the expectation is taken with respect to the distribution $(\mb P^\star)^\infty$ governing the sample path and the empirical estimator. The mean squared error penalizes the mismatch between the actual cost $c(x, \mb P^\star)$ and its estimator $\hat c(x,\hat {\mb P}_T)$. Events in which we are left disappointed ($c(x, \mb P^\star) > \hat c(x,\hat {\mb P}_T)$) are not treated differently from positive surprises ($c(x, {\mb P^\star}) < \hat c(x,\hat {\mb P}_T)$). In a decision-making context where the goal is to minimize costs, however, disappointments (underestimated costs) are more harmful than positive surprises (overestimated costs). While statisticians strive for accuracy by minimizing a symmetric estimation error, decision makers endeavor to limit the one-sided prediction disappointment. 

\begin{definition}[Out-of-sample disappointment]
  \label{def:oos-disappointment}
  For any data-driven predictor $\hat c$ the probability 
  \begin{subequations}
    \label{eq:oos-disappointment}
    \begin{equation}
      \label{eq:prediction-disappointment}
      \mb P^\infty \left( c(x, \mb P) > \hat c(x, \hat{\mb P}_T)\right)
    \end{equation}
    is referred to as the out-of-sample prediction disappointment of~$x\in X$ under model $\mb P\in\mc P$. Similarly, for any data-driven prescriptor $\hat x$ induced by a data-driven predictor $\hat c$ the probability 
    \begin{equation}
      \label{eq:prescription-disappointment}
      \mb P^\infty \left( c(\hat x(\hat{\mb P}_T), \mb P) > \hat c(\hat x(\hat {\mb P}_T), \hat {\mb P}_T)\right)
    \end{equation}
  \end{subequations}
  is termed the out-of-sample prescription disappointment under model $\mb P\in\mc P$.
\end{definition}

The out-of-sample prediction disappointment quantifies the probability (with respect to the sample path distribution $\mb P^\infty$ under some model $\mb P\in\mc P$) that the expected cost $c(x, \mb P)$ of a fixed decision $x$ exceeds the predicted cost $\hat c(x, \hat{\mb P}_T)$. Thus, the out-of-sample prediction disappointment is independent of the actual realization of the empirical estimator $\hat {\mb P}_T$ but depends on the hypothesized model $\mb P$. A similar statement holds for the out-of-sample prescription disappointment.

The main objective of this paper is to construct attractive data-driven predictors and prescriptors, which are optimal in a sense to be made precise below. We first develop a notion of optimality for data-driven predictors and extend it later to data-driven prescriptors. As indicated above, a crucial requirement for any data-driven predictor is that it must limit the out-of-sample disappointment. This informal requirement can be operationalized either in an asymptotic sense or in a finite sample sense.
\begin{itemize}
\item[(i)] {\bf Asymptotic guarantee:} 
  As $T$ grows, the out-of-sample prediction disappointment~\eqref{eq:prediction-disappointment} decays exponentially at a rate at least equal to $r> 0$ up to first order in the exponent, that is,
  \begin{equation}
    \label{eq:asymptotic-prediction-guarantee}
    \limsup_{T\to\infty} \frac{1}{T} \log \mb P^\infty \left( c(x, \mb P) > \hat c(x, \hat{\mb P}_T)\right) \leq -r \quad \forall x\in X,\;\mb P\in\mc P.
  \end{equation}
\item[(ii)] {\bf Finite sample guarantee:} For every fixed $T$, the out-of-sample prediction disappointment~\eqref{eq:prediction-disappointment} is bounded above by a {\em known} function $g(T)$ that decays exponentially at rate at least equal to $r>0$ to first order in the exponent, that is,
  \begin{equation}
    \label{eq:finite-sample-prediction-guarantee}
    \mb P^\infty \left( c(x, \mb P) > \hat c(x, \hat{\mb P}_T)\right) \leq g(T)~~ \forall x\in X,\;\mb P\in\mc P, \; T\in\mb N,
  \end{equation}
  where $\limsup_{T\to\infty} \frac{1}{T} \log g(T) \leq -r$.
\end{itemize}
The inequalities \eqref{eq:asymptotic-prediction-guarantee} and \eqref{eq:finite-sample-prediction-guarantee} are imposed across all models $\mb P\in\mc P$. This ensures that they are satisfied under the true model $\mb P^\star$, which is only known to reside within $\mc P$. By requiring the inequalities to hold for all $x\in X$, we further ensure that the out-of-sample prediction disappointment is eventually small irrespective of the chosen decision. Note that the finite sample guarantee~\eqref{eq:finite-sample-prediction-guarantee} is sufficient but not necessary for the asymptotic guarantee \eqref{eq:asymptotic-prediction-guarantee}. Knowing the finite sample bounds $g(T)$ has the advantage, amongst others, that one can determine the sample complexity
\[
  \min\left\{T_0\in\mb N: g(T)\leq \beta,~ \forall T\geq T_0\right\},
\]
that is, the minimum number of samples needed to certify that the out-of-sample prediction disappointment does not exceed a prescribed significance level $\beta\in [0,1]$.

At first sight the requirements~\eqref{eq:asymptotic-prediction-guarantee} and \eqref{eq:finite-sample-prediction-guarantee} may seem restrictive, and the existence of data-driven predictors with exponentially decaying out-of-sample disappointment may be questioned. Below we will argue, however, that these requirements are in fact natural and satisfied by all reasonable predictors. To see this, note that if the training data is generated by $\mb P$, then the empirical distribution $\hat{\mb P}_T$ converges $\mb P^\infty$-almost surely to $\mb P$ by virtue of the strong law of large numbers. Thus, the out-of-sample disappointment of a predictor $\hat c$ with $\hat c(x,{\mb P})>c(x,\mb P)$ must decay to 0 as $T$ grows. Conversely, if $\hat c(x,{\mb P})<c(x,\mb P)$, then the out-of-sample disappointment of $\hat c$ must approach~1 as $T$ tends to infinity. The following example shows that the out-of-sample disappointment generically fails to vanish asymptotically in the limiting case when $\hat c(x,{\mb P})=c(x,\mb P)$.

\begin{example}[Large out-of-sample disappointment]
  \label{ex:large-disappointment}
  Set the cost function to $\gamma(x,\xi)=\xi$. In this case, the sample average predictor approximates the expected cost $c(x,\mb P)= \sum_{i\in\Xi} i\mb P(i)$ by its sample mean $c(x,\hat{\mb P}_T)=\frac{1}{T}\sum_{t=1}^T \xi_t$. As the sample size $T$ tends to infinity, the central limit theorem implies that \[\sqrt{T} [c(x,\hat{\mb P}_T)-c(x,\mb P)]\] converges in law to a normal distribution with mean $0$ and variance $\mb E_{\mb P}[(\xi - \mb E_{\mb P}[\xi])^2]$. Thus,
  \[
    \lim_{T\to\infty}\mb P^\infty \left( c(x, \mb P) > \hat c(x, \hat{\mb P}_T)\right) =\lim_{T\to\infty}\mb P^\infty \left( \sqrt{T}\left(\hat c(x, \hat{\mb P}_T)-c(x,\mb P)\right)<0\right) = \frac{1}{2},
  \]
  which means that the out-of-sample prediction disappointment remains large for all sample sizes. The sample average predictor hence violates the asymptotic guarantee~\eqref{eq:asymptotic-prediction-guarantee} and the stronger finite sample guarantee~\eqref{eq:finite-sample-prediction-guarantee}. Note that by adding {\em any} positive constant to the sample average predictor, we recover a predictor with exponentially decaying out-of-sample disappointment.
\end{example}

In the following we call a predictor $\hat c$ {\em conservative} if $\hat c(x,{\mb P'})>c(x,\mb P')$ for all decisions $x\in X$ and estimator realizations $\mb P'\in\mc P$. The above discussion shows that if we require the out-of-sample disappointment to decay asymptotically, we must focus on conservative predictors. Basic results from large deviations theory further ensure that the out-of-sample disappointment of any conservative predictor necessarily decays at an exponential rate. Specifically, asymptotic guarantees of the type~\eqref{eq:asymptotic-prediction-guarantee} hold whenever the empirical distribution $\hat {\mb P}_T$ satisfies a {\em weak large deviation principle}, while finite sample guarantees of the type~\eqref{eq:finite-sample-prediction-guarantee} hold when $\hat {\mb P}_T$ satisfies a {\em strong large deviation principle}. As will be shown in Section \ref{sec:ldt}, the empirical distribution does satisfy weak and strong large deviation principles. One predictor that fails to be conservative is the sample average predictor.

For ease of exposition, we henceforth denote by $\mc C$ the set of all data-driven predictors, that is, all continuous functions that map $X\times\mc P$ to the reals. Moreover, we introduce a partial order $\preceq_\mc C$ on $\mc C$ defined through
\[
  \hat c_1\preceq_\mc C \hat c_2 \quad \iff\quad \hat c_1(x,\mb P')\leq \hat c_2(x,\mb P')\quad\forall x\in X,\, \mb P'\in\mc P
\]
for any $\hat c_1,\hat c_2\in\mc C$. Thus, $\hat c_1\preceq_\mc C \hat c_2$ means that $\hat c_1$ is (weakly) less conservative than $\hat c_2$. The problem of finding the least conservative predictor among all data-driven predictors whose out-of-sample disappointment decays at rate at least $r> 0$ can thus be formalized as the following {\em vector optimization problem}.
\begin{equation}
  \label{eq:optimal-predictor}
  \begin{array}{ll}
    \displaystyle\mathop{\text{minimize}}_{\hat c\in\mc C}{}_{\preceq_\mc C} & \hat c \\
    \text{subject to} & \displaystyle \limsup_{T\to\infty} \frac{1}{T} \log \mb P^\infty \left( c(x, \mb P) > \hat c(x, \hat{\mb P}_T)\right)\leq -r \quad \forall x\in X,\;\mb P\in\mc P
  \end{array}
\end{equation}
We highlight that the minimization in~\eqref{eq:optimal-predictor} is understood with respect to the partial order $\preceq_\mc C$. Thus, the relation $\hat c_1\preceq_\mc C \hat c_2$ between two feasible decision means that $\hat c_1$ is weakly preferred to $\hat c_2$. However, not all pairs of feasible decisions are comparable, that is, it is possible that both $\hat c_1\npreceq_\mc C \hat c_2$ and $\hat c_2\npreceq_\mc C \hat c_1$. A predictor~$\hat c{}^\star$ is a {\em strongly} optimal solution for~\eqref{eq:optimal-predictor} if it is feasible and weakly preferred to every other feasible solution ({\em i.e.}, every $\hat c\neq \hat c{}^\star$ feasible in~\eqref{eq:optimal-predictor} satisfies $\hat c{}^\star\preceq_\mc C \hat c$). Similarly, $\hat c{}^\star$ is a {\em weakly} optimal solution for~\eqref{eq:optimal-predictor} if it is feasible and if every other solution preferred to $\hat c{}^\star$ is infeasible ({\em i.e.}, every $\hat c\neq \hat c{}^\star$ with $\hat c\preceq_\mc C \hat c{}^\star$ is infeasible in~\eqref{eq:optimal-predictor}). {While vector optimization problems can have many weak solutions, we point out that strong solutions are necessarily unique. To see this, assume for the sake of contradiction that $\hat c_1^\star$ and $\hat c_2^\star$ are two strong solutions of~\eqref{eq:optimal-predictor}. In this case the strong optimality of $\hat c_1^\star$ implies that $\hat c_2^\star \preceq_\mc C \hat c_1^\star$, while the strong optimality of $\hat c_2^\star$ implies that $\hat c_1^\star \preceq_\mc C \hat c_2^\star$. These two relations imply that $\hat c_1^\star = \hat c_2^\star$, that is, there cannot be two different strongly optimal solutions.}

We are now ready to construct a meta-optimization problem akin to~\eqref{eq:optimal-predictor}, which enables us to identify the best prescriptor. To this end, we henceforth denote by~$\mc X$ the set of all data-driven predictor-prescriptor-pairs $(\hat c,\hat x)$, where $\hat c\in\mc C$, and $\hat x$ is a prescriptor induced by $\hat c$ as per Definition~\ref{def:dd_prediction}. Moreover, we equip $\mc X$ with a partial order $\preceq_\mc X$, which is defined through
\[
  (\hat c_1,\hat x_1) \preceq_\mc X (\hat c_2,\hat x_2) \quad \iff\quad \hat c_1(\hat x_1(\mb P'),\mb P')\leq \hat c_2(\hat x_2(\mb P') ,\mb P')\quad\forall \mb P'\in\mc P.
\]
Note that $\hat c_1\preceq_\mc C \hat c_2$ actually implies $(\hat c_1,\hat x_1) \preceq_\mc X (\hat c_2,\hat x_2)$ but not vice versa. The problem of finding the least conservative predictor-prescriptor-pair whose out-of-sample prescription disappointment decays at rate at least $r> 0$ can now be formalized as the following vector optimization problem.
\begin{equation}
  \label{eq:optimal-prescriptor}
  \begin{array}{ll}
    \displaystyle\mathop{\text{minimize}}_{(\hat c,\hat x)\in\mc X}{}_{\preceq_\mc X} & (\hat c,\hat x) \\
    \text{subject to} & \displaystyle \limsup_{T\to\infty} \frac{1}{T} \log \mb P^\infty \left( c(\hat x(\hat{\mb P}_T), \mb P) > \hat c(\hat x(\hat {\mb P}_T), \hat {\mb P}_T)\right) \leq -r \quad \forall \mb P\in\mc P
  \end{array}
\end{equation}
Generic vector optimization problems typically only admit weak solutions. In Section \ref{sec:distr_rob_optimization} we will show, however, that \eqref{eq:optimal-predictor} as well as \eqref{eq:optimal-prescriptor} admit (unique) strong solutions in closed form. In fact, we will show that these closed-form solutions have a natural interpretation as the solutions of convex distributionally robust optimization problems.

{
\begin{remark}[Out-of-sample and in-sample performance]
The natural performance measure to quantify the goodness of a data-driven prescriptor $\hat x$ is its {\em out-of-sample performance} $c(\hat x(\hat{\mb P}_T), \mb P^\star)$ under the true model~$\mb P^\star$. As $\mb P^\star$ is unknown, however, the out-of-sample performance cannot be optimized directly. A na\"ive remedy would be to formulate a meta-optimization problem that minimizes the worst-case (or some average) of the out-of-sample performance of $\hat x$ across all models $\mb P\in \mc P$. The approach proposed here optimizes the out-of-sample performance implicitly. Indeed, the meta-optimization problem~\eqref{eq:optimal-prescriptor} represents $\hat x$ as a minimizer of some predictor $\hat c$, where $\hat c(\hat x(\hat {\mb P}_T), \hat {\mb P}_T)$ should be interpreted as the {\em in-sample performance} of $\hat x$. Instead of minimizing the out-of-sample performance of $\hat x$, problem~\eqref{eq:optimal-prescriptor} minimizes the in-sample performance of $\hat x$ but ensures through the constraints on the disappointment that the out-of-sample performance is smaller than the in-sample performance with increasingly high confidence as the sample size grows. In this sense, problem~\eqref{eq:optimal-prescriptor} minimizes a tight upper bound on the out-of-sample performance of $\hat x$. 
\end{remark}}

\section{Large deviation principles}
\label{sec:ldt}

Large deviations theory provides bounds on the exact exponential rate at which the probabilities of atypical estimator realizations decay under a model~$\mb P$ as the sample size $T$ tends to infinity. These bounds are expressed in terms of the relative entropy of $\hat{\mb P}_T$ with respect to~$\mb P$.

\begin{definition}[Relative entropy] 
  \label{def:relative_entropy}
  The relative entropy of an estimator realization $\mb P'\in\mc P$ with respect to a model $\mb P\in\mc P$ is defined as
  \begin{align*}
    \D{\mb P'}{\mb P} = \sum_{i\in\Xi} \mb P'(i) \log\left(\frac{\mb P'(i)}{\mb P(i)}\right),
  \end{align*}
  where we use the conventions $0 \log(0/p)=0$ for any $p\geq 0$ and $p' \log(p'/0)=\infty$ for any $p'> 0$.
\end{definition}

The relative entropy is also known as information for discrimination, cross-entropy, information gain or Kullback-Leibler divergence \citep{kullback1951information}. The following proposition summarizes the key properties of the relative entropy relevant for this paper.

\begin{proposition}[Relative entropy]
  \label{prop:relative-entropy}
  The relative entropy enjoys the following properties:
  \begin{itemize}
    \item[(i)] {\bf Information inequality:} $\D{\mb P'}{\mb P}\geq 0$ for all $\mb P,\mb P'\!\in\mc P$, while $\D{\mb P'}{\mb P}= 0$ if and only if $\mb P'=\mb P$.
  \item[(ii)] {\bf Convexity:} For all pairs $(\mb P'_1, \mb P_1), (\mb P'_2, \mb P_2)\in \mc P\times \mc P$ and $\lambda\in[0,1]$ we have
    \[
      \D{(1-\lambda)\mb P'_1 +\lambda \mb P'_2}{(1-\lambda)\mb P_1 +\lambda \mb P_2}\leq (1-\lambda) \D{\mb P'_1}{\mb P_1}+\lambda \D{\mb P'_2}{\mb P_2}.
    \] 
  \item[(iii)] {\bf Lower semicontinuity} $\D{\mb P'}{\mb P}\geq 0$ is lower semicontinuous in $ (\mb P', \mb P)\in \mc P\times \mc P$.
  \end{itemize}
  
\end{proposition}
\proof{Proof.} 
 Assertions~(i) and~(ii) follow from Theorems~2.6.3 and~2.7.2 in \citet{cover2006elements}, respectively, while assertion~(iii) follows directly from the definition of the relative entropy and our standard conventions regarding the natural logarithm. \hfill $\square$
\endproof

We now show that the empirical estimators satisfy a weak {\em large deviation principle} (LDP). This result follows immediately from a finite version of Sanov's classical theorem. A textbook proof using the so-called method of types can be found in~\citet[Theorem~11.4.1]{cover2006elements}. As the proof is illuminating and to keep this paper self-contained, we sketch the proof in Appendix~\ref{sec:proofs}. 

\begin{theorem}[Weak LDP]
  \label{thm:sldp-iid}
  If the samples $\{\xi_t\}_{t\in\mb N}$ are drawn independently from some~$\mb P\in\mc P$, then for every Borel set $\mc D\subseteq \mc P$ the sequence of empirical distributions $\{\hat {\mb P}_T\}_{T\in\mc P}$ satisfies
  \begin{subequations}
    \label{eq:ldp_exponential_rates}
     \begin{equation}
     \label{eq:ldp_exponential_rates_ub}
        \limsup_{T\to \infty}~\frac1T \log \mb P^\infty( \hat {\mb P}_T \in \mc D ) \leq -\inf_{\mb P' \in \mc D} \, \D{\mb P'}{\mb P}.
    \end{equation}
     If additionally $\mb P>0$, then for every Borel set $\mc D\subseteq \mc P$ we have\footnote{Here, the interior of $\mc D$ is taken with respect to the subspace topology on $\mathcal P$. Recall that a set $\mc D\subseteq \mc P$ is open in the subspace topology on $\mc P$ if $\mc D=\mc P\cap \mc O$ for some set $\mc O\subseteq \Re^d$ that is open in the Euclidean topology on $\Re^d$.}
     \begin{equation}
      \label{eq:ldp_exponential_rates_lb}
      \liminf_{T\to \infty}~\,\frac1T \log \mb P^{\infty}( \hat {\mb P}_T \in \mc D )\geq -\inf_{\mb P' \in \interior \mc D} \, \D{\mb P'}{\mb P}.
     \end{equation}
  \end{subequations}
\end{theorem}

Note that the inequality~\eqref{eq:ldp_exponential_rates_ub} provides an {\em upper} LDP bound on the exponential rate at which the probability of the event $\hat {\mb P}_T \in \mc D$ decays under model $\mb P$. This upper bound is expressed in terms of a convex optimization problem that minimizes the relative entropy of $\mb P'$ with respect to $\mb P$ across all estimator realizations $\mb P'$ within $\mc D$. Similarly,~\eqref{eq:ldp_exponential_rates_lb} offers a {\em lower} LDP bound on the decay rate. Note that in~\eqref{eq:ldp_exponential_rates_lb} the relative entropy is minimized over the {\em interior} of~$\mc D$ instead of $\mc D$.

\begin{figure}
   \centering
   \includegraphics[width=0.45\textwidth]{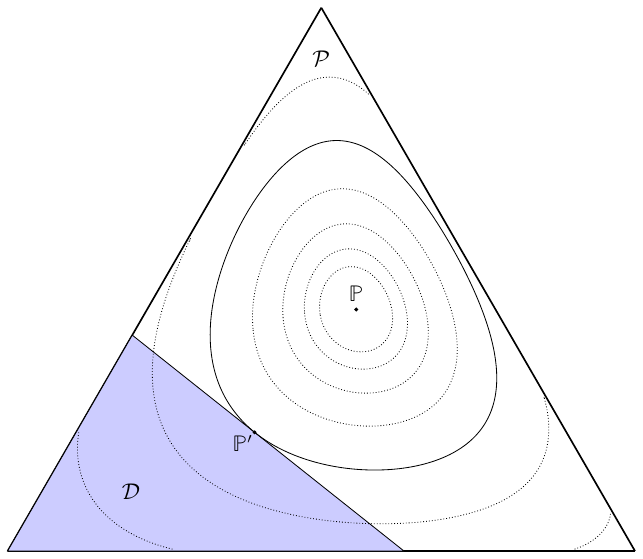}
   %\subimport{figs/}{LDTsimplex.tex}
   \caption{Visualization of the LDP~\eqref{eq:ldp_exponential_rates}. If $\mc D\subseteq \mc P$ is $\rm I$-continuous and $\mb P\notin\mc D$, then the probability $\mb P^\infty(\hat {\mb P}_T\in \mc D)$ decays at the exponential rate $\inf_{\mb P' \in \mc D} \D{\mb P'}{\mb P}$, which can be viewed as the relative entropy distance of $\mb P$ from $\mc D$.}
   \label{fig:ldp}
 \end{figure}

If the data-generating model $\mb P$ itself belongs to~$\mc D$, then $\inf_{\mb P' \in \mc D}  \D{\mb P'}{\mb P}=\D{\mb P}{\mb P}=0$, which leads to the trivial upper bound $\mb P^\infty(\hat {\mb P}_T \in \mc D )\leq 1$. On the other hand, if $\mc D$ has empty interior ({\em e.g.}, if $\mc D=\{\mb P\}$ is a singleton containing only the true model), then $\inf_{\mb P' \in \interior{\mc D}} I(\mb P', \mb P)=\infty$, which leads to the trivial lower bound $\mb P^\infty(\hat {\mb P}_T \in \mc D )\geq 0$. Non-trivial bounds are obtained if $\mb P\notin \mc D$ and $\interior\mc D\neq \emptyset$. In these cases the relative entropy bounds the exponential rate at which the probability of the atypical event $\hat {\mb P}_T\in\mc D$ decays with $T$. For some sets $\mc D$ this rate of decay is precisely determined by the relative entropy. Specifically, a Borel set $\mc D\subseteq \mc P$ is called $\rm I$-continuous under model $\mb P$ if
\begin{equation*}
  \inf_{\mb P' \in \interior{\mc D}} \, \D{\mb P'}{\mb P}= \inf_{\mb P' \in \mc D} \,  \D{\mb P'}{\mb P}.
\end{equation*}
Clearly, every open set $\mc D\subseteq \mc P$ is $\rm I$-continuous under any model $\mb P$. Moreover, as the relative entropy is continuous in $\mb P'$ for any fixed $\mb P>0$, every Borel set $\mc D\subseteq \mc P$ with $\mc D \subseteq \cl(\interior(\mc D))$ is $\rm I$-continuous under $\mb P$ whenever $\mb P>0$. The LDP~\eqref{eq:ldp_exponential_rates} implies that for large $T$ the probability of an $\rm I$-continuous set $\mc D$ decays at rate $\inf_{\mb P'\in \mc D} \D{\mb P'}{\mb P}$ under model $\mb P$ to first order in the exponent, that is, we have
\begin{equation}
  \label{eq:ldp:continuous}
  \mb P^\infty( \hat {\mb P}_T \in \mc D ) = e^{-T \inf_{\mb P'\in \mc D} \D{\mb P'}{\mb P}+o(T)}.
\end{equation}
If we interpret the relative entropy $\D{\mb P'}{\mb P}$ as the distance of $\mb P$ from $\mb P'$, then the decay rate of $\mb P^\infty( \hat {\mb P}_T \in \mc D )$ coincides with the distance of the model $\mb P$ from the atypical event set $\mc D$; see Figure~\ref{fig:ldp}.  Moreover, if $\mc D$ is $\rm I$-continuous under $\mb P$, then \eqref{eq:ldp:continuous} implies that $\mb P^\infty(\hat {\mb P}_T \in \mc D)\leq \beta$ whenever
\[
  T\gtrsim \frac{1}{r} \cdot \log \left(\frac{1}{\beta}\right),
\]
where $r = \inf_{\mb P'\in\mc D}\, \D{\mb P'}{\mb P}$ is the $\rm I$-distance from $\mb P$ to the set $\mc D$, and $\beta\in(0,1)$ is a prescribed significance level.

The weak LDP of Theorem~\ref{thm:sldp-iid} provides only {\em asymptotic} bounds on the decay rates of atypical events. However, one can also establish a {\em strong} LDP, which offers {\em finite sample guarantees}.  Most results of this paper, however, are based on the weak LDP of Theorem~\ref{thm:sldp-iid}.

\begin{theorem}[Strong LDP]
  \label{thm:sldp-iid-strong}
   If the samples $\{\xi_t\}_{t\in\mb N}$ are drawn independently from some~$\mb P\in\mc P$, then for every Borel set $\mc D\subseteq \mc P$ the sequence of empirical distributions $\{\hat {\mb P}_T\}_{T\in\mc P}$ satisfies
  \begin{align}
    \label{eq:sldp}
    \mb P^\infty\left( \hat {\mb P}_T \in \mc D \right) & \leq (T+1)^{d} e^{-T \inf_{\mb P' \in \mc D} \,  \D{\mb P'}{\mb P}} \quad \forall T \in \mb N.
  \end{align} 
\end{theorem}
\proof{Proof.}
  The claim follows immediately from inequality~\eqref{eq:sanov:ub} in the proof of Theorem~\ref{thm:sldp-iid} in Appendix~\ref{sec:proofs}. Note that~\eqref{eq:sanov:ub} does not rely on the assumption that $\mb P>0$. \hfill $\square$
\endproof

\section{Distributionally robust predictors and prescriptors are optimal}
\label{sec:distr_rob_optimization}

Armed with the fundamental results of large deviations theory, we now endeavor to identify the least conservative data-driven predictors and prescriptors whose out-of-sample disappointment decays at a rate no less than some prescribed threshold $r> 0$ under any model $\mb P\in \mc P$, that is, we aim to solve the vector optimization problems \eqref{eq:optimal-predictor} and \eqref{eq:optimal-prescriptor}. 

\subsection{Distributionally robust predictors}
\label{ssec:dd_prediction}
The relative entropy lends itself to constructing a data-driven predictor in the sense of Definition~\ref{def:dd_prediction}. We will show below that this predictor is strongly optimal in~\eqref{eq:optimal-predictor}.

\begin{definition}[Distributionally robust predictors]
  \label{def:dro-predictor}
  For any fixed threshold $r\geq 0$, we define the data-driven predictor $\hat c_{r}:X\times\mc P\to \Re$ through
  \begin{equation}
    \label{eq:gen:dro}
    \hat c_{r}(x, \mb P') = \sup_{\mb P \in \mc P}\, \left\{ c(x, \mb P) : \D{\mb P'}{ \mb P}\leq r \right\} \qquad \forall x\in X,~\mb P'\in\mc P.
  \end{equation}
\end{definition}

The data-driven predictor $\hat c_r$ admits a distributionally robust interpretation. In fact, $\hat c_r(x, \mb P')$ represents the worst-case expected cost associated with the decision $x$, where the worst case is taken across all models $\mb P\in\mc P$ whose relative entropy distance to $\mb P'$ is at most $r$. Observe that the supremum in~\eqref{eq:gen:dro} is always attained because $c(x,\mb P)$ is linear in $\mb P$ and the feasible set of~\eqref{eq:gen:dro} is compact, which follows from the compactness of $\mc P$ and the lower semicontinuity of the relative entropy in $\mb P$ for any fixed $\mb P'$; see Proposition~\ref{prop:relative-entropy}(iii). Note also that $\hat c_r(x, \mb P')$ can be evaluated efficiently because \eqref{eq:gen:dro} constitutes a convex conic optimization problem with $d$ decision variables. A particularly simple and efficient method to evaluate $\hat c_r(x, \mb P')$ is to solve the one-dimensional convex minimization problem dual to~\eqref{eq:gen:dro} by using bisection or another line search method.

\begin{proposition}[Dual representation of $\hat c_r$]
  \label{thm:strong-dual-final-finite}
  If $r>0$ and $\bar\gamma(x)=\max_{i\in \Xi} \, \gamma(x, i)$ denotes the worst-case cost function, then the distributionally robust predictor admits the dual representation
  \begin{equation}
    \label{eq:strong-dual-final-finite}
      \hat c_{r}(x,\mb P')= \min_{\alpha\ge \bar\gamma(x)} \alpha - e^{-r} \prod_{i\in \Xi} \left(\alpha-\gamma(x, i) \right)^{\mb P'(i)}.
    \end{equation}
  Problem~\eqref{eq:strong-dual-final-finite} has a minimizer $\alpha^\star$ that satisfies $\bar\gamma(x) \leq \alpha^\star \leq \frac{\bar\gamma(x)-e^{-r} c(x, \mb P')}{1-e^{-r}}$.
\end{proposition}
\proof{Proof.}
See Appendix \ref{sec:proofs}. \hfill $\square$
\endproof

\begin{remark}[Sample average predictor]
  For $r=0$ the distributionally robust predictor $\hat c_r$ collapses to the sample average predictor of Example~\ref{ex:naive-predictor}. Indeed, because of the strict positivity of the relative entropy $\D{\mb P'}{\mb P} >0$ for $\mb P'\neq \mb P$, see Proposition~\ref{prop:relative-entropy}(i), we have that 
  \[
    \hat c_0(x, \mb P') = c(x, \mb P').
  \]
  As shown in Example~\ref{ex:large-disappointment}, the sample average predictor fails to offer asymptotic or finite sample guarantees of the form~\eqref{eq:asymptotic-prediction-guarantee} and~\eqref{eq:finite-sample-prediction-guarantee}, respectively. 
\end{remark}

\begin{remark}[Alternative distributionally robust predictors]
  The relative entropy can also be used to construct a reverse distributionally robust predictor $\check c_r\in\mc C$ defined through
  \begin{equation}
    \label{eq:KL-wrong}
    \check c_r(x, \mb P') = \sup_{\mb P \in \mc P}\, \left\{ c(x, \mb P) : \D{\mb P}{\mb P'}\leq r \right\}\qquad \forall x\in X,~\mb P'\in\mc P.
  \end{equation}
  In contrast to $\hat c_r$, the reverse distributionally robust predictor $\check c_r$ fixes the {\em second} argument of the relative entropy and maximizes over the {\em first} argument. Note that $\check c_r$ can be viewed as the entropic value-at-risk of the uncertain cost $\gamma(x,\xi)$; see \citep[Theorem~3.3]{ahmadi2012entropic}. Another predictor related to $\hat c$ is the restricted distributionally robust predictor $\bar c_r\in\mc C$ defined through
    \begin{equation}
      \label{eq:KL-restricted}
      \bar c_r(x, \mb P') = \sup_{\mb P \in \mc P}\, \left\{ c(x, \mb P) : \mb P \ll \mb P' ,~  \D{\mb P'}{\mb P}\leq r \right\}\qquad \forall x\in X,~\mb P'\in\mc P,
    \end{equation}
    where $\mb P \ll \mb P'$ expresses the requirement that $\mb P$ must be absolutely continuous with respect to $\mb P'$. Formally, $\mb P \ll \mb P'$ means that $\mb P(i)=0$ for all outcomes $i \in \Xi$ with $\mb P'(i)=0$. By \citep[Definition~5.1]{ahmadi2012entropic}, $\bar c_r$ can be interpreted as the negative log-entropic risk of $\gamma(x,\xi)$.
  
  The predictors $\hat c_r$ and $\check c_r$ differ because the relative entropy fails to be symmetric. We emphasize that the reverse predictor $\check c_r$ has appeared often in the literature on distributionally robust optimization, see, {\em e.g.}, \citep{ben2013robust, calafiore2007ambiguous, hu2013kullback, lam2016robust, wang2009likelihood}. The predictors $\hat c_r$ and $\bar c_r$ differ, too, because of the additional constraint $\mb P \ll \mb P'$, which is significant when not all outcomes in $\Xi$ have been observed.  The statistical properties of the predictor $\bar c_r$ have been analyzed by \citet{lam2016} and more recently by \citet{duchi2016statistics} from the perspective of the empirical likelihood theory introduced by \citet{owen1988empirical}.  The predictor $\hat c_r$ suggested here has not yet been studied extensively even though---as we will demonstrate below---it displays attractive theoretical properties that are not shared by~either $\check c_r$ or $\bar c_r$. The difference between $\hat c_r$ and $\check c_r$ or $\bar c_r$ is significant. Indeed, both $\check c_r$ and $\bar c_r$ hedge only against models $\mb P$ that are absolutely continuous with respect to the (observed realization of the) empirical distribution~$\mb P'$. While it is clear that the empirical distribution must be absolutely continuous with respect to the data-generating distribution, however, the converse implication is generally false. Indeed, an outcome can have positive probability even if it does not show up in a given finite time series. By taking the worst case only over models that are absolutely continuous with respect to $\mb P'$, both predictors $\check c_r$ and $\bar c_r$ potentially ignore many models that could have generated the observed data.
\end{remark}

We first establish that $\hat c_r$ indeed belongs to the set $\mc C$ of all data-driven predictors, that is, the family of continuous functions mapping $X\times \mc P$ to the reals. 

\begin{proposition}[Continuity of $\hat c_r$]
  \label{prop:continuity_cr}
  If $r\ge 0$, then the distributionally robust predictor~$\hat c_r$ is continuous on $X\times\mc P$.
\end{proposition}
\proof{Proof.}
By Proposition~\ref{thm:strong-dual-final-finite}, the distributionally robust predictor $\hat c_r$ admits the dual representation~\eqref{eq:strong-dual-final-finite}. Note that the objective function of~\eqref{eq:strong-dual-final-finite} is manifestly continuous in $(\alpha,x,\mb P')$ and that~\eqref{eq:strong-dual-final-finite} is guaranteed to have a minimizer in the compact interval $[\bar\gamma(x), \frac{\bar\gamma(x)-e^{-r} c(x, \mb P')}{1-e^{-r}}]$, whose boundaries depend continuously on $(x,\mb P')$. Consequently, the predictor $\hat c_r$ is continuous by Berge's celebrated maximum theorem \citep[pp.~115--116]{berge1963topological}. \hfill $\square$
\endproof

We now analyze the performance of the distributionally robust data-driven predictor $\hat c_r$ using arguments from large deviations theory. 
The parameter $r$ encoding the predictor $\hat c_{r}$ captures the fundamental trade-off between out-of-sample disappointment and accuracy, which is inherent to any approach to data-driven prediction. Indeed, as~$r$ increases, the predictor $\hat c_r$ becomes more reliable in the sense that its out-of-sample disappointment decreases. However, increasing~$r$ also results in more conservative (pessimistically biased) predictions. In the following we will demonstrate that $\hat c_r$ strikes indeed an optimal balance between reliability and conservatism. 

\begin{theorem}[Feasibility of $\hat c_r$]
  \label{thm:pred:feasibility}
  If $r\geq 0$, then the predictor $\hat c_r$ is feasible in~\eqref{eq:optimal-predictor}. 
\end{theorem}
\proof{Proof.}
  From Proposition~\ref{prop:continuity_cr} we already know that $\hat c_r\in\mc C$. It remains to be shown that the out-of-sample disappointment of $\hat c_r$ decays at a rate of at least $r$.
  We have $c(x, \mb P) > \hat c_{r} (x, \hat {\mb P}_T)$ if and only if the estimator $\hat {\mb P}_T$ falls within the disappointment set
  \[
    \mc D(x, \mb P)= \left\{ \mb P'\in \mc P: c(x, \mb P) >  \hat c_{r}(x, \mb P') \right\}.
  \]
  Note that by the definition of $\hat c_{r}$, we have 
  \[
    \D{\mb P'}{\mb P}\leq r\quad \implies \quad \hat c_{r}(x,\mb P')=\sup_{\mb P''\in\mc P} \left\{ c(x,\mb P''): I(\mb P',\mb P'')\leq r\right\}\geq c(x,\mb P). 
  \]
  By contraposition, the above implication is equivalent to
  \[
    c(x,\mb P)> \hat c_r(x,\mb P')\quad \implies \quad \D{\mb P'}{\mb P}> r.
  \]
  Therefore, $\mc D(x, \mb P)$ is a subset of
  \[
    \mc I(\mb P)= \left \{ \mb P'\in \mc P:  \D{\mb P'}{\mb P}>r \right\}
  \]
  irrespective of $x\in X$. We thus have
  \begin{align*} 
    \limsup_{T\to\infty}\frac 1T \log \mb P^\infty \left( \hat {\mb P}_T \in {\mc D}(x, \mb P) \right) \leq\limsup_{T\to\infty} \frac 1T \log \mb P^\infty \left( \hat {\mb P}_T \in \mc I(\mb P) \right)\leq  -\inf_{\mb P' \in \mc I(\mb P)} \, \D{\mb P'}{\mb P}  \leq-r, 
  \end{align*}
  where the first inequality holds because ${\mc D}(x, \mb P)\subseteq \mc I(\mb P)$, while the second inequality exploits the weak LDP upper bound \eqref{eq:ldp_exponential_rates_ub}. Thus, $\hat c_r$ is feasible in \eqref{eq:optimal-predictor}. \hfill$\square$
\endproof

{Note that any predictor $\hat c$ with $\hat c_r \preceq_{\mathcal C} \hat c$ has a smaller disappointment set than $\hat c_r$, and thus the out-of-sample disappointment of $\hat c$ decays at least as fast as that of $\hat c_r$. Hence, $\hat c$ is also feasible in~\eqref{eq:optimal-predictor}. In particular, this immediately implies that if we inflate the relative entropy ball of the distributionally robust predictor $\hat c_r$ to any larger ambiguity set, we obtain another predictor that is feasible in~\eqref{eq:optimal-predictor}. As an example, consider the total variation predictor
  \begin{equation*}
    \hat c^{\rm{tv}}_r(x, \mb P') = \sup_{\mb P \in \mc P}\, \left\{ c(x, \mb P) :  \norm{\mb P- \mb P'}_{\rm{tv}}\leq \sqrt{2r} \right\}\qquad \forall x\in X,~\mb P'\in\mc P,
  \end{equation*}
  where $\norm{\mb P-\mb P'}_{\rm{tv}}$ denotes the total variation distance between $\mb P$ and $\mb P'$. Pinsker's classical inequality asserts that $\norm{\mb P-\mb P'}_{\rm{tv}}\leq \sqrt{2 \,\D{\mb P'}{\mb P}}$ for all $\mb P$ and $\mb P'$ in $\mc P$. Thus, we have $\hat c_r \preceq_{\mathcal C} \hat c^{\rm{tv}}_r$, which implies that the total variation predictor is feasible in \eqref{eq:optimal-predictor}. This suggests that~\eqref{eq:optimal-predictor} has a rich feasible set.
}

The following main theorem establishes that $\hat c_r$ is not only a feasible but even a strongly optimal solution for the vector optimization problem~\eqref{eq:optimal-predictor}. This means that if an arbitrary data-driven predictor $\hat c$ predicts a lower expected cost than $\hat c_r$ even for a single estimator realization $\mb P'\in\mc P$, then $\hat c$ must suffer from a higher out-of-sample disappointment than $\hat c_r$ to first order in the exponent.

\begin{theorem}[Optimality of $\hat c_r$]
  \label{thm:optimality}
  If $r>0$, then $\hat c_r$ is strongly optimal in~\eqref{eq:optimal-predictor}. 
\end{theorem}
\proof{Proof.}
  Assume for the sake of argument that $\hat c_r$ fails to be a strong solution for~\eqref{eq:optimal-predictor}. Thus, there exists a data-driven predictor $\hat c\in\mc C$ that is feasible in~\eqref{eq:optimal-predictor} but not dominated by $\hat c_r$, that is, $\hat c_r\npreceq_\mc C \hat c$. This means that there exists $x\in X$ and $\mb P'_0\in\mc P$ with $\hat c_r(x,\mb P'_0)>\hat c(x,\mb P'_0)$. For later reference we set $\epsilon=\hat c_r(x,\mb P'_0)-\hat c(x,\mb P'_0)>0$. In the remainder of the proof we will demonstrate that $\hat c$ cannot be feasible in \eqref{eq:optimal-predictor}, which contradicts our initial assumption.

  Let $\mb P_0\in\mc P$ be an optimal solution of problem~\eqref{eq:gen:dro} at $\mb P'=\mb P'_0$. Thus, we have $I(\mb P'_0,\mb P_0)\leq r$ and
  \begin{equation}
    \label{eq:inequality2}
    \hat c_r(x, \mb P'_0)= c(x, \mb P_0).
  \end{equation}
  In the following we will first perturb $\mb P_0$ to obtain a model $\mb P_1$ that is $\frac{\epsilon}{2}$-suboptimal in~\eqref{eq:gen:dro} but satisfies $I(\mb P'_0,\mb P_1)< r$. Subsequently, we will perturb $\mb P_1$ to obtain a model $\mb P_2$ that is $\epsilon$-suboptimal in~\eqref{eq:gen:dro} but satisfies $I(\mb P'_0,\mb P_2)< r$ as well as $\mb P_2>0$.
  
  To construct $\mb P_1$, consider all models $\mb P(\lambda)=\lambda\mb P'_0 + (1-\lambda)\mb P_0$, $\lambda\in [0,1]$, on the line segment between $\mb P'_0$ and $\mb P_0$. As $r$ is strictly positive, the convexity of the relative entropy implies that
  \[
  	\D{\mb P_0'}{\mb P(\lambda)}\leq  \lambda \D{\mb P_0'}{\mb P_0'}+ (1-\lambda)\D{\mb P_0'}{\mb P_0}\leq(1-\lambda) r< r \quad \forall \lambda\in(0,1].
  \]
  Moreover, as the expected cost $c(x,\mb P(\lambda))$ changes continuously
 in $\lambda$, there exists a sufficiently small $\lambda_1\in (0,1]$ such that $\mb P_1=\mb P(\lambda_1)$ and $r_1=\D{\mb P'_0}{\mb P_1}$ satisfy $0<r_1<r$ and
  \begin{equation*}
    c(x, \mb P_0)< c(x, \mb P_1) +\frac{\epsilon}{2}.
  \end{equation*}
  To construct $\mb P_2$, we consider all models $\mb P(\lambda)=\lambda\mb U + (1-\lambda) \mb P_1$, $\lambda\in [0,1]$, on the line segment between the uniform distribution $\mb U$ on $\Xi$ and $\mb P_1$. By the convexity of the relative entropy we have
  \[
  	\D{\mb P_0'}{\mb P(\lambda)}\leq  \lambda \D{\mb P_0'}{\mb U}+ (1-\lambda)\D{\mb P_0'}{\mb P_1}\leq  \lambda \D{\mb P_0'}{\mb U}+ (1-\lambda)r_1\quad \forall \lambda\in[0,1].
  \]
  As $r_1<r$ and the expected cost $c(x,\mb P(\lambda))$ changes continuously in $\lambda$, there exists a sufficiently small $\lambda_2\in (0,1]$ such that $\mb P_2=\mb P(\lambda_2)$ and $r_2=\D{\mb P'_0}{\mb P_2}$ satisfy $0<r_2<r$, $\mb P_2>0$ and
  \begin{equation}
    \label{eq:inequality3}
    c(x, \mb P_0)< c(x, \mb P_2) +\epsilon.
  \end{equation}
  In summary, we thus have
  \begin{align}
    \label{eq:theta'-feasibility}
    \hat c(x,\mb P'_0) = \hat c_r(x,\mb P'_0)-\epsilon =  c(x,\mb P_0) -\epsilon < c(x,\mb P_2)\leq \hat c_r(x,\mb P'_0),
  \end{align}
  where the first equality follows from the definition of $\epsilon$, and the second equality exploits~\eqref{eq:inequality2}. Moreover, the strict inequality holds due to~\eqref{eq:inequality3}, and the weak inequality follows from the definition of $\hat c_r$ and the fact that $\D{\mb P'_0}{\mb P_2}=r_2< r$.

  In the remainder of the proof we will argue that the prediction disappointment $\mb P_{2}^\infty (c(x, \mb P_2) > \hat c(x,\hat {\mb P}_T))$ under model $\mb P_2$ decays at a rate of at most $r_2<r$ as the sample size~$T$ tends to infinity. In analogy to the proof of Theorem~\ref{thm:pred:feasibility}, we define the set of disappointing estimator realizations as
  \[
    \mc D(x,\mb P_2) = \set{\mb P' \in \mc P}{c(x, \mb P_2) > \hat c(x, \mb P')}.
  \]
  This set contains $\mb P'_0$ due to the strict inequality in~\eqref{eq:theta'-feasibility}. Moreover, as $\hat c\in\mc C$ is continuous, $ \mc D(x,\mb P_2)$ is an open subset of $\mc P$.  Thus, we find
  \[
    \inf_{\mb P'\in \interior \mc D(x,\mb P_2)} \, \D{\mb P'}{\mb P_2} = \inf_{\mb P'\in \mc D(x,\mb P_2)} \, \D{\mb P'}{\mb P_2} \leq \D{\mb P'_0}{\mb P_2}= r_2,
  \]
  where the inequality holds because $\mb P'_0\in \mc D(x,\mb P_0)$, and the last equality follows from the definition of $r_2$. As the empirical distributions $\{\hat{\mb P}_T\}_{T\in\mb N}$ obey the LDP lower bound \eqref{eq:ldp_exponential_rates_lb} under $\mb P_2>0$, we finally conclude that 
  \[
    -r<-r_2\leq-\inf_{\mb P'\in \interior \mc D(x,\mb P_2)} \, \D{\mb P'}{\mb P_2} \leq \liminf_{T\rightarrow \infty}\frac 1T \log \mb P_{2}^\infty \left(  \hat {\mb P}_T \in \mc D(x,\mb P_2) \right). \]
  The above chain of inequalities implies, however, that $\hat c$ is infeasible in problem~\eqref{eq:optimal-predictor}. This contradicts our initial assumption, and thus, $\hat c_r$ must indeed be a strong solution of \eqref{eq:optimal-predictor}. \hfill $\square$
\endproof

Theorem \ref{thm:optimality} asserts that the distributionally robust predictor $\hat c_r$ is optimal among all data-driven predictors representable as continuous functions of the empirical distribution $\hat {\mb P}_T$. That is, any attempt to make it less conservative invariably increases the out-of-sample prediction disappointment. We remark that the class of predictors which depend on the data only through $\hat{\mb P}_T$ is vast. These predictors constitute arbitrary continuous functions of the data that are independent of the order in which the samples were observed. As the samples are independent and identically distributed, there are in fact no meaningful data-driven predictors that display a more general dependence on the data. 

Note that in the above discussion all guarantees are fundamentally asymptotic in nature. Using Theorem~\ref{thm:sldp-iid-strong} one can show, however, that $\hat c_r$ also satisfies finite sample guarantees. 

\begin{theorem}[Finite sample guarantee]
  \label{thm:pred:finite_sample}
  The out-of-sample disappointment of the distributionally robust predictor $\hat c_r$ enjoys the following finite sample guarantee under any model $\mb P\in\mc P$ and for any~$x\in X$. 
  \begin{equation}
    \label{eq:pred:finite_sample}
    \mb P^\infty \left( c(x, \mb P) > \hat c_{r} (x, \hat {\mb P}_T) \right) \leq  (T+1)^d e^{-rT} \quad \forall T\in\mb N
  \end{equation}
\end{theorem}
\proof{Proof.}
The proof of this result widely parallels that of Theorem~\ref{thm:pred:feasibility} but uses the strong LDP upper bound~\eqref{eq:sldp} in lieu of the weak upper bound~\eqref{eq:ldp_exponential_rates_ub}. Details are omitted for brevity.  \hfill $\square$
\endproof

\subsection{Distributionally robust prescriptors}
\label{ssec:dd_prescription}

The distributionally robust predictor $\hat c_r$ of Definition \ref{def:dro-predictor} induces a corresponding prescriptor.

\begin{definition}[Distributionally robust prescriptors]
  \label{def:dro-prescriptor}
  Denote by $\hat c_r$, $r\geq 0$, the distributionally robust data-driven predictor of Definition \ref{def:dro-predictor}. We can then define the data-driven prescriptor $\hat x_r:\mc P\to X$ as a quasi-continuous function satisfying
  \begin{equation}
    \label{eq:gen:dro_prescriptor}
    \hat x_r(\mb P') \in \arg\min_{x\in X} \, \hat c_r(x, \mb P') \quad \forall \mb P'\in\mc P.
  \end{equation}
\end{definition}

Note that the minimum in~\eqref{eq:gen:dro_prescriptor} is attained because $X$ is compact and $\hat c_r$ is continuous due to Proposition~\ref{prop:continuity_cr}. Thus, there exists at least one function $\hat x_r$ satisfying~\eqref{eq:gen:dro_prescriptor}. In the next proposition we argue that this function can be chosen to be quasi-continuous as desired. 

\begin{proposition}[Quasi-continuity of $\hat x_r$]
  \label{prop:measurability_xr}
  If $r\geq 0$, then there exists a quasi-continuous data-driven predictor $\hat x_r$ satisfying \eqref{eq:gen:dro_prescriptor}.
\end{proposition}
\proof{Proof.}
    Denote by $\Gamma(\mb P')= \arg\min_{x\in X} \, \hat c_r(x, \mb P')$ the argmin-mapping of problem~\eqref{eq:gen:dro}, and observe that $\Gamma(\mb P')$ is compact and non-empty for every $\mb P'\in\mc P$ because $\hat c_r$ is continuous and $X$ is compact. As $X$ is independent of $\mathbb P'$, Berge's maximum theorem \citep[pp.~115--116]{berge1963topological} further implies that $\Gamma$ is upper semicontinuous. As $\mc P$ is a Baire space and $X$ is a metric space, \citep[Corollary~4]{Matejdes:87} finally guarantees that there exists a quasi-continuous function $\hat x_r:\mc P\to X$ with $\hat x_r(\mb P')\in\Gamma(\mb P')$ for all $\mb P'\in\mc P$. \hfill $\square$
\endproof

Propositions~\ref{prop:continuity_cr} and~\ref{prop:measurability_xr} imply that $(\hat c_r,\hat x_r)$ belongs to the family $\mc X$ of all data-driven predictor-prescriptor-pairs. Using a similar reasoning as in Theorem~\ref{thm:pred:feasibility}, we now demonstrate that the out-of-sample disappointment of $\hat x_r$ decays at rate at least $r$ as $T$ tends to infinity. Thus, $\hat x_r$ provides trustworthy prescriptions.

\begin{theorem}[Feasibility of $(\hat c_r, \hat x_r)$]
  \label{thm:presc:feasibility}
  If $r\geq 0$, then the predictor-prescriptor-pair $(\hat c_r,\hat x_r)$ is feasible in~\eqref{eq:optimal-prescriptor}. 
\end{theorem}

\proof{Proof.}
  Propositions~\ref{prop:continuity_cr} and~\ref{prop:measurability_xr} imply that $(\hat c_r,\hat x_r)\in\mc X$. It remains to be shown that the out-of-sample disappointment of $\hat x_r$ decays at a rate of at least $r$. To this end, define $\mc D(x,\mb P)$ and $\mc I(\mb P)$ as in the proof of Theorem \ref{thm:pred:feasibility}, and recall that $\mc D(x,\mb P)\subseteq \mc I(\mb P)$ for every decision $x\in X$ and model $\mb P\in\mc P$. Thus, for every fixed estimator realization $\mb P'\in\mc P$ the following implication holds
  \begin{align*}
    c(\hat x_{r}(\mb P'), \mb P) > \hat c_{r} (\hat x_{r}(\mb P'), \mb P')\quad \implies\quad &\exists x\in X \text{ with } c(x, \mb P) > \hat c_{r} (x, \mb P')\\
    \implies \quad & \mb P' \in \cup_{x\in X} \mc D(x,\mb P) \\ \implies \quad & \mb P' \in \mc I(\mb P),
  \end{align*}
  which in turn implies
  \[
    \limsup_{T\to\infty} \frac 1T \log \mb P^\infty \left( c(\hat x_{r}(\hat {\mb P}_T), \mb P) > \hat c_{r} (\hat x_{r}(\hat {\mb P}_T), \hat {\mb P}_T) \right) \leq \limsup_{T\to\infty}\frac 1T \log \mb P^\infty \left( \hat {\mb P}_T \in \mc I(\mb P) \right) \leq -r
  \]
  for every model $\mb P\in\mc P$. Note that the second inequality in the above expression has already been established in the proof of Theorem~\ref{thm:pred:feasibility}. Thus, the claim follows.  \hfill $\square$

Next, we argue that $(\hat c_r,\hat x_r)$ is a strongly optimal solution for the vector optimization problem~\eqref{eq:optimal-prescriptor}.

\begin{theorem}[Optimality of $(\hat c_r, \hat x_r)$]
  \label{thm:optimality_prescriptor}
  If $r>0$, then $(\hat c_r,\hat x_r)$ is strongly optimal in~\eqref{eq:optimal-prescriptor}. 
\end{theorem}
\proof{Proof.}
Assume for the sake of argument that $(\hat c_r, \hat x_r)$ fails to be a strong solution for~\eqref{eq:optimal-prescriptor}. Thus, there exists a data-driven prescriptor $(\hat c, \hat x)\in\mc X$ that is feasible in~\eqref{eq:optimal-prescriptor} but not dominated by $(\hat c_r, \hat x_r)$, that is, $(\hat c_r, \hat x_r)\npreceq_{\mc X} (\hat c, \hat x)$. This means that there exists $\mb P_0'\in\mc P$ with $\hat c_r(\hat x_r(\mb P'_0), \mb P'_0)>\hat c(\hat x(\mb P'_0),\mb P'_0)$. As $X$ is compact and $\hat c$ is continuous, the cost $\hat c(\hat x(\mb P'),\mb P')$ of the prescriptor $\hat x$ under the corresponding predictor $\hat c$ is continuous in~$\mb P'$ \citep[pp.~115--116]{berge1963topological}. Similarly, $\hat c_r(\hat x_r(\mb P'),\mb P')$ is continuous in~$\mb P'$. Recall also that $\hat x$ is quasi-continuous and therefore continuous on a dense subset of $\mc P$ \citep{Bledsoe1952}. Thus, we may assume without loss of generality that $\hat x$ is continuous at $\mb P_0'$. For later reference we set $\epsilon=\hat c_r(\hat x(\mb P'_0),\mb P'_0)-\hat c(\hat x(\mb P'_0),\mb P'_0)>0$.

In the remainder of the proof we will demonstrate that $(\hat c, \hat x)$ cannot be feasible in \eqref{eq:optimal-prescriptor}, which contradicts our initial assumption. To this end, let $\mb P_0\in\mc P$ be an optimal solution of problem~\eqref{eq:gen:dro} at $x=\hat x(\mb P'_0)$ and $\mb P'=\mb P'_0$. Thus, we have $I(\mb P'_0,\mb P_0)\leq r$ and
  \begin{equation}
    \label{eq:prescr:inequality2}
    \hat c_r(\hat x(\mb P'_0), \mb P'_0)= c(\hat x(\mb P'_0), \mb P_0).
  \end{equation}
  Next, we first perturb $\mb P_0$ to obtain a model $\mb P_1$ that is strictly $\frac{\epsilon}{2}$-suboptimal in~\eqref{eq:gen:dro} but satisfies $I(\mb P'_0,\mb P_1)=r_1< r$. Subsequently, we perturb $\mb P_1$ to obtain a model $\mb P_2$ that is strictly $\epsilon$-suboptimal in~\eqref{eq:gen:dro} but satisfies $I(\mb P'_0,\mb P_2)=r_2< r$ as well and $\mb P_2>0$. The distributions $\mb P_1$ and $\mb P_2$ can be constructed exactly as in the proof of Theorem~\ref{thm:optimality}. Details are omitted for brevity. Thus, we find 
  \begin{align}
    \label{eq:prescr:theta'-feasibility}
    \hat c(\hat x(\mb P'_0),\mb P'_0) = \hat c_r(\hat x(\mb P'_0),\mb P'_0)-\epsilon =  c(\hat x(\mb P'_0),\mb P_0) -\epsilon < c(\hat x(\mb P'_0),\mb P_2)\leq \hat c_r(\hat x(\mb P'_0),\mb P'_0),
  \end{align}
  where the first equality follows from the definition of $\epsilon$, and the second equality exploits~\eqref{eq:prescr:inequality2}. Moreover, the strict inequality holds because $\mb P_2$ is strictly $\epsilon$-suboptimal in~\eqref{eq:gen:dro}, while the weak inequality follows from the definition of $\hat c_r$ and the fact that $I(\mb P'_0,\mb P_2)=r_2< r$.

  It remains to be shown that the prediction disappointment $\mb P_{2}^\infty (c(\hat x(\hat{\mb P}_T), \mb P_2) > \hat c(\hat x(\hat{\mb P}_T),\hat {\mb P}_T))$ under model $\mb P_2$ decays at a rate of at most $r_2<r$ as the sample size~$T$ tends to infinity. To this end, we define the set of disappointing estimator realizations as
  \[
    \mc D(\mb P_2) = \set{\mb P' \in \mc P}{c(\hat x(\mb P'), \mb P_2) > \hat c(\hat x(\mb P'), \mb P')}.
  \]
  This set contains $\mb P'_0$ due to the strict inequality in~\eqref{eq:prescr:theta'-feasibility}. Recall now that $\hat x$ is continuous at $\mb P'=\mb P'_0$ due to our choice of $\mb P'_0$. As the predictors $\hat c$ and $\hat c_r$ are both continuous on their entire domain, the compositions $\hat c(\hat x(\mb P'), \mb P')$ and $c(\hat x(\mb P'), \mb P_2)$ are both continuous at $\mb P'=\mb P'_0$. This implies that $\mb P'_0$ belongs actually to the interior of $\mc D(\mb P_2)$. Thus, we find
  \[
    \inf_{\mb P'\in \interior \mc D(\mb P_2)} \, \D{\mb P'}{\mb P_2}\leq \D{\mb P'_0}{\mb P_2}= r_2,
  \]
  where the last equality follows from the definition of $r_2$. As the empirical distributions $\{\hat{\mb P}_T\}_{T\in\mb N}$ obey the LDP lower bound \eqref{eq:ldp_exponential_rates_lb} under $\mb P_2>0$, we finally conclude that 
  \[
    -r<-r_2\leq-\inf_{\mb P'\in \interior \mc D(\mb P_2)} \, \D{\mb P'}{\mb P_2} \leq \liminf_{T\rightarrow \infty}\frac 1T \log \mb P_{2}^\infty \left(  \hat {\mb P}_T \in \mc D(\mb P_2) \right). \]
  The above chain of inequalities implies, however, that $(\hat c, \hat x)$ is infeasible in problem~\eqref{eq:optimal-prescriptor}. This contradicts our initial assumption, and thus, $(\hat c_r, \hat x_r)$ must indeed be a strong solution of \eqref{eq:optimal-prescriptor}. \hfill $\square$
\endproof

All guarantees discussed so far are asymptotic in nature. As in the case of the predictor $\hat c_r$, however, the prescriptor $\hat x_r$ can also be shown to satisfy finite sample guarantees. 

\begin{theorem}[Finite sample guarantee]
  \label{thm:presc:finite_sample}
  The out-of-sample disappointment of the distributionally robust prescriptor $\hat x_{r}$ enjoys the following finite sample guarantee under any model $\mb P\in\mc P$. 
  \begin{equation}
    \label{eq:prescription:guarantee}
    \mb P^\infty \left( c(\hat x_r(\hat {\mb P}_T), \mb P) > \hat c_r (\hat x_{r}(\hat {\mb P}_T), \hat {\mb P}_T) \right) \leq (T+1)^{d} e^{-rT} \quad \forall T\in\mb N
  \end{equation}
\end{theorem}
\proof{Proof.}
The proof of this result parallels those of Theorems~\ref{thm:pred:feasibility} and~\ref{thm:presc:feasibility} but uses the strong LDP upper bound~\eqref{eq:sldp} in lieu of the weak upper bound~\eqref{eq:ldp_exponential_rates_ub}. Details are omitted for brevity.  \hfill $\square$
\endproof

We stress that the finite sample guarantees of Theorems~\ref{thm:pred:finite_sample} and \ref{thm:presc:finite_sample} as well as the strong optimality properties portrayed in Theorems~\ref{thm:optimality} and~\ref{thm:optimality_prescriptor} are independent of a particular dataset. They guarantee that $\hat c_r$ and $\hat x_r$ provide trustworthy predictions and prescriptions, respectively, {\em before} the data is revealed.

  \begin{remark}[Optimal hypothesis testing]
    \citet{bertsimas2018data} propose to construct predictors and prescriptors from statistical hypothesis tests. A hypothesis test uses i.i.d.\ samples $\xi_1,\dots, \xi_T$ drawn from the unknown true distribution $\mb P^\star$ to decide whether the null hypothesis $\mb P^\star=\mb P$ is false for a fixed model $\mb P\in\mc P$. Specifically, the null hypothesis is rejected (it is declared that $\mb P^\star\neq \mb P$) if the empirical distribution $\hat{\mb P}_T$ associated with the observed sample path falls outside of a (measurable) acceptance region $A_T(\mb P)\subseteq \mc P$, which depends on the conjectured model $\mb P$ and the sample size $T$. Otherwise, it is deemed that there is insufficient data to reject the null hypothesis. 
    
    \citet{bertsimas2018data} associate with each hypothesis test a predictor
    \begin{equation}
    \label{eq:hypothesis-predictor}
      \hat c(x, \mb P') = \left\{\begin{array}{r@{\hspace{0.75em}}l}
        \sup &  c(x, \mb P )\\[0.5em]
        \st & \mb P \in \mc P\\[0.5em]
                 & \mb P' \in A_T(\mb P),
      \end{array}\right.
    \end{equation}
    which evaluates the worst-case expected cost across all models $\mb P\in\mc P$ that pass the hypothesis test in view of the realization $\mb P'\in\mc P_T=\mc P\cap \{0, 1/T, \dots, (T-1)/T,1\}^d$ of the empirical distribution~$\hat {\mb P}_T$.
    
    The quality of a hypothesis test is usually measured by its type~I error $\mb P^\infty(\hat{\mb P}_T \not\in A_T(\mb P))$, that is, the probability of falsely rejecting the null hypothesis, as well as its type II error $\mb Q^\infty(\hat{\mb P}_T \in A_T(\mb P))$, that is, the probability of falsely accepting the null hypothesis if the data follows a distribution $\mb Q\neq \mb P$. A particularly popular test is the likelihood ratio test, which uses the acceptance region
    \[
      A^\star_T(\mb P)= \set{\mb P'\in\mc P_T}{\textstyle {\mb P^\infty(\hat{\mb P}_T = \mb P')}/{\sup_{\mb Q\neq \mb P} \mb Q^\infty(\hat{\mb P}_T = \mb P')}\geq e^{-rT}}.
    \]
    \citet{zeitouni1992generalized} prove that the likelihood ratio test is optimal in the following sense. Among all hypothesis tests whose type~I error decays at a rate of at least $r$, $\limsup_{T\to\infty} \frac 1T \log \mb P^\infty(\hat{\mb P}_T \notin A_T(\mb P)) \leq -r$, the likelihood ratio test minimizes the nagative decay rate of the type~II error $\limsup_{T\to\infty} \frac 1T \log \mb Q^\infty(\hat{\mb P}_T \in A_T(\mb P))$ simultaneously for all models $\mb Q\in \mc P$ with $\mb Q\neq \mb P$. \citet[Theorem 11.1.2]{cover2006elements} further establish that the likelihood ratio of an estimator realization $\mb P'\in\mc P_T$ under two alternative distributions $\mb Q$ and $\mb P$ satisfies $\log({\mb P^\infty(\hat{\mb P}_T = \mb P')}/{\mb Q^\infty(\hat{\mb P}_T = \mb P')}) = -T ( \D{\mb P'}{\mb P} - \D{\mb P'}{\mb Q})$. The acceptance region of the likelihood ratio test thus simplifies to
    \[
      A_T^\star(\mb P) =  \set{\mb P'\in\mc P_T}{\textstyle \D{\mb P'}{\mb P} \leq r + \inf_{\mb Q\neq \mb P} \D{\mb P'}{\mb Q}}
      =  \set{\mb P'\in\mc P_T}{ \D{\mb P'}{\mb P} \leq r},
    \]
    where the equality holds because $\inf_{\mb Q\neq \mb P} \D{\mb P'}{\mb Q}=0$. Hence, the distributionally robust predictor~$\hat c_r$ that is strongly optimal in the meta-optimization problem~\eqref{eq:optimal-predictor} coincides with the hypothesis test-based predictor~\eqref{eq:hypothesis-predictor} corresponding to the likelihood ratio test. 
\end{remark}

\section{Extension to continuous state spaces}
\label{sec:continuous}
Assume now that the realizations of the random parameter $\xi$ may range over an arbitrary compact set $\Xi\subseteq \Re^d$ that is not necessarily finite. In analogy to the discrete case, we denote by $\mc P$ the family of all Borel probability distributions supported on $\Xi$. Note that $\mc P$ is now a convex subset of an infinite-dimensional space, which significantly complicates the problem of finding optimal predictors and prescriptors. We equip $\mc P$ with the standard topology of weak convergence of distributions, recalling that the weak topology is metrized by the Prokhorov metric \citep{prokhorov_convergence_1956}. Consequently, we equip $X\times \mc P$ with the product of the standard Euclidean topology on $X$ and the weak topology on $\mc P$. In the remainder of this section we analyze to what extent---and under what additional conditions---the results for finite state spaces carry over to the more general continuous case. As this analysis requires more subtle mathematical techniques, we relegate all proofs to Appendix~\ref{sec:proofs}.

We first note that the definitions of model-based predictors and prescriptors require no changes. In order to evaluate the expectation in the definition of the model-based predictor $c(x,\mb P)=\int_\Xi \gamma(x,\xi) \d \mb P(\xi)$, however, we now need to evaluate an integral with respect to $\mb P$ instead of a finite sum. Throughout this section we assume that the cost function $\gamma(x,\xi)$ is jointly continuous in $x$ and $\xi$. This implies via the compactness of $X$ and $\Xi$ that $c(x, \mb P)$ is continuous in $x$ and $\mb P$, which in turn guarantees that a model-based prescriptor $x^\star(\mb P)\in\arg\min_{x\in X}c(x,\mb P)$ exists for every $\mb P\in\mc P$.

\begin{lemma}[Continuity of model-based predictors]
  \label{lemma:r=0}
  If $\gamma(x,\xi)$ is continuous on the compact set $X\times \Xi$, then $c(x, \mb P)$ is continuous on $X\times \mc P$.
\end{lemma}

As in the case of a discrete state space, we study data-driven predictors and prescriptors that depend on the training data $\{\xi_t\}_{t=1}^T$ only through the empirical distribution. Because $\Xi$ may now have infinite cardinality, we redefine the empirical distribution as $\hat {\mb P}_T= \frac{1}{T} \sum_{t=1}^T \delta_{\xi_t}$, where~$\delta_{\xi_t}$ denotes the Dirac point mass at $\xi_t$. Using this new  definition of $\hat {\mb P}_T$, we then define data-driven predictors and prescriptors exactly as in Section~\ref{ssec:dd_predictions_prescriptions}. As $\Xi$ is compact, one can show that $\mc P$ is compact in the weak topology \citep{prokhorov_convergence_1956}. Moreover, as the weak topology is metrized by the Prokhorov metric, $\mc P$ constitutes a (locally) compact metric space. The Baire category theorem thus implies that $\mc P$ is a Baire space \citep{baire1899fonctions}. Corollary~4 in~\citep[p.~120]{Matejdes:87}, which applies because $\mc P$ is a Baire space and $X$ is a metric space, further ensures that for any valid (continuous) predictor $\hat c$ the set-valued mapping $\arg\min_{x\in X}\, \hat c(x, {\mb P'})$ admits a quasi-continuous selector $\hat x$, which serves as a valid data-driven prescriptor. Using the exact same reasoning as in Section~\ref{ssec:dd_predictions_prescriptions}, one can show that the points of continuity of any quasi-continuous prescriptor are dense in $\mc P$.

The best predictors and predictor-prescriptor-pairs can again be found by solving the meta-optimization problems~\eqref{eq:optimal-predictor} and~\eqref{eq:optimal-prescriptor}, respectively. In order to construct near-optimal solutions for these meta-optimization problems, we recall the definition of the relative entropy between arbitrary distributions $\mb P'$ and $\mb P$ on a compact set $\Xi\subseteq \Re^d$. 

\begin{definition}[Generalized relative entropy] 
  \label{def:generalized_relative_entropy}
  The relative entropy of $\mb P'\in\mc P$ with respect to $\mb P\in\mc P$ is defined as
  \begin{align*}
    \D{\mb P'}{\mb P} =
    \begin{cases}
      \int_{\Xi} \log\left({\d\mb P'}/{\d\mb P}(\xi)\right) \d \mb P'(\xi) & {\rm{if}~} \mb P' \ll \mb P,\\[0.5em]
      +\infty & {\rm{otherwise,}}
    \end{cases}
  \end{align*}
  where $\mb P' \ll \mb P$ means that $\mb P'$ is absolutely continuous with respect to $\mb P$, while ${\d\mb P'}/{\d\mb P}(\xi)$ denotes the Radon-Nikodym derivative of $\mb P'$ with respect to $\mb P$, which exists if $\mb P'\ll\mb P$ \citep{nikodym_sur_1930}.
\end{definition}

The properties of the relative entropy portrayed in Proposition~\ref{prop:relative-entropy} hold verbatim in the more general setting considered here \citep{vanerven2014renyi}. Using the generalized definition of the relative entropy, the distributionally robust predictor $\hat c_r$ and the corresponding prescriptor $\hat x_r$ can be constructed as in Definitions~\ref{def:dro-predictor} and~\ref{def:dro-prescriptor}, respectively. In the following we will show that the predictor $\hat c_r$ is continuous, which ensures that the prescriptor $\hat x_r$ can always be chosen to be quasi-continuous. To this end, we first derive a dual representation for $\hat c_r$.

\begin{proposition}[Dual representation revisited]
  \label{thm:strong-dual-final}
  If $r>0$ and $\bar\gamma(x)=\max_{\xi\in \Xi} \, \gamma(x, \xi)$ is the worst-case cost function, then the distributionally robust predictor admits the dual representation
  \begin{equation}
    \label{eq:strong-dual}
      \hat c_{r}(x,\mb P')=  \min_{\alpha\ge\bar\gamma(x)} \alpha - e^{-r}\cdot \exp\left( \int_{\Xi} \log\left(\alpha-\gamma(x, \xi)\right) \d\mb P'(\xi)\right).
  \end{equation}
 Problem~\eqref{eq:strong-dual} has a minimizer $\alpha^\star \leq \frac{\bar\gamma(x) -e^{-r} c(x, \mb P')}{1-e^{-r}}$.
\end{proposition}

Proposition~\ref{thm:strong-dual-final} extends Proposition~\ref{thm:strong-dual-final-finite} to compact continuous state spaces and suggests that~$\hat c_{r}(x,\mb P')$ can be computed via bisection or other line search methods. Thus, the computational tractability of problem~\eqref{eq:strong-dual} largely hinges on our ability to efficiently evaluate the geometric mean $\exp\left( \int_{\Xi} \log\left(\alpha-\gamma(x, \xi)\right) \d\mb P'(\xi)\right)$ for any fixed $\alpha$. For example, if $\mb P'$ coincides with (a realization of) the empirical distribution $\hat {\mb P}_T$, we recover the geometric mean of $\alpha-\gamma(x, \xi)$ along a sample path, { which can be reformulated as the optimal value of a tractable second-order cone program involving $\mc O(T)$ constraints and auxiliary variables \citep[Section~6.2.3.5]{nesterov1994interior}.}
\[
  \exp\left( \int_{\Xi} \log\left(\alpha-\gamma(x, \xi)\right) \d\mb P'(\xi)\right) = \left(\prod_{t=1}^T \left(\alpha-\gamma(x, \xi_t) \right) \right)^{1/T}
\]
To our best knowledge, the dual representation~\eqref{eq:strong-dual} is new. The closest result we are aware of is the dual representation of the negative log-entropic risk measure derived in \citep[Theorem~5.1]{ahmadi2012entropic}. Indeed, the negative log-entropic risk of $\gamma(x,\xi)$ coincides with the restricted distributionally robust predictor $\bar c_{r}(x, \mb P')$. Recall from~\eqref{eq:KL-restricted} that $\bar c_{r}(x, \mb P')$ differs from $c_{r}(x, \mb P')$ only in that it imposes the additional constraint $\mb P\ll \mb P'$ when evaluating the worst-case expected cost. Using Theorem~5.1 by \citet{ahmadi2012entropic} one can thus show that the dual representation of $\bar c_{r}(x, \mb P')$ differs from~\eqref{eq:strong-dual} only in that it replaces $\bar \gamma(x)$ with $\inf\{\bar\gamma: \mb P'[ \gamma(x,\xi)\le\bar \gamma]=1\} \le\bar\gamma(x)$. Maybe surprisingly, however, the derivation of~\eqref{eq:strong-dual} provided here is substantially more challenging.

\begin{proposition}[Continuity of $\hat c_r$ revisited]
  \label{prop:pred:feasibility:revisited}
  If $r\ge 0$, then the distributionally robust predictor~$\hat c_r$ is continuous on $X\times\mc P$.
\end{proposition}

Proposition~\ref{prop:pred:feasibility:revisited} ensures that $\hat c_r\in\mc C$. As any continuous predictor induces a quasi-continuous prescriptor, we may thus conclude that there exists a valid distributionally robust prescriptor~$\hat x_r$ such that~$(\hat c_r, \hat x_r)\in\mc X$. It now only remains to establish that these predictors and predictor-prescriptor-pairs are the unique strong solutions of the meta-optimization problems~\eqref{eq:optimal-predictor} and~\eqref{eq:optimal-prescriptor}, respectively. In Section~\ref{sec:distr_rob_optimization} this was achieved by leveraging the weak LDP portrayed in Theorem~\ref{thm:sldp-iid}. Luckily, this LDP carries over to the more general setting considered here---albeit with a subtle difference.

\begin{theorem}[Weak LDP revisited]
  \label{thm:sldp-iid-continous}
  If the samples $\{\xi_t\}_{t\in\mb N}$ are drawn independently from some~$\mb P\in\mc P$, then for every set $\mc D\subseteq \mc P$ the sequence of empirical distributions $\{\hat {\mb P}_T\}_{T\in\mb N}$ satisfies
  \begin{subequations}
    \label{eq:ldp_exponential_rates-continuous}
     \begin{align}
     \label{eq:ldp_exponential_rates_ub-continuous}
        \limsup_{T\to \infty}~\frac1T \log \mb P^\infty( \hat {\mb P}_T \in \mc D )& \leq -\inf_{\mb P' \in \cl \mc D} \, \D{\mb P'}{\mb P},\\
      \label{eq:ldp_exponential_rates_lb-continuous}
      \liminf_{T\to \infty}~\,\frac1T \log \mb P^{\infty}( \hat {\mb P}_T \in \mc D )& \geq -\inf_{\mb P' \in \interior \mc D} \, \D{\mb P'}{\mb P}.
     \end{align}
  \end{subequations}
\end{theorem}
\proof{Proof.}
See~\citet[Section~2]{Csiszar2006}.
\hfill $\square$
\endproof

Formally, Theorem~\ref{thm:sldp-iid-continous} is almost identical to Theorem~\ref{thm:sldp-iid}. However, the weak LDP upper bound~\eqref{eq:ldp_exponential_rates_ub-continuous} differs from~\eqref{eq:ldp_exponential_rates_ub} in that the minimization over all estimator realizations $\mb P'$ on the right hand side runs over the {\em closure} of $\mc D$. This subtle difference invalidates the proof of Theorem~\ref{thm:pred:feasibility}, and thus we need a new approach to show that $\hat c_r$ is feasible in~\eqref{eq:optimal-predictor}. Moreover, the weak LDP lower bound~\eqref{eq:ldp_exponential_rates_lb-continuous} does {\em not} rely on any structural assumptions about $\mb P$. Note that the condition $\mb P>0$ in Theorem~\ref{thm:sldp-iid} was only imposed for convenience to simplify the proof of~\eqref{eq:ldp_exponential_rates_lb} in Appendix~\ref{sec:proofs}.

As in the case of finite state spaces, one can now show that the distributionally robust predictor~$\hat c_r$ is the unique strong solution of the meta-optimization problem~\eqref{eq:optimal-predictor}.

\begin{theorem}[Feasibility and optimality of $\hat c_r$ revisited]
  \label{thm:pred:feasibility-revisited}
  If $r\geq 0$, then the predictor $\hat c_r$ is feasible in~\eqref{eq:optimal-predictor}. Moreover, if $r>0$, then $\hat c_r$ is strongly optimal in~\eqref{eq:optimal-predictor}.
\end{theorem}

While we did not manage to prove that $(\hat c_r,\,\hat x_r)$ is feasible in the meta-optimization problem~\eqref{eq:optimal-prescriptor}, we still could show that it is {\em essentially} feasible and strongly optimal in a precise sense. 

  \begin{theorem}[Feasibility and optimality of $(\hat c_r, \hat x_r)$ revisited]
  \label{thm:presc:feasibility-revisited}
  If $r\geq 0$, then the shifted predictor-prescriptor-pair $(\hat c_r+\epsilon,\hat x_r)$ is feasible in~\eqref{eq:optimal-prescriptor} for every $\epsilon>0$. Moreover, if $r>0$, then $(\hat c_r, \hat x_r)$ is preferred to every feasible solution of~\eqref{eq:optimal-prescriptor}---even though it may be infeasible.
\end{theorem}

Theorem~\ref{thm:presc:feasibility-revisited} asserts that $(\hat c_r, \hat x_r)$ is less conservative than any predictor-prescriptor pair feasible in the meta-optimization problem~\eqref{eq:optimal-prescriptor} and that $(\hat c_r, \hat x_r)$ can be made feasible in~\eqref{eq:optimal-prescriptor} by shifting the distributionally robust predictor $\hat c_r$ up by just a tiny amount. For practical purposes this means that $(\hat c_r, \hat x_r)$ is indeed essentially optimal. Whether $(\hat c_r,\,\hat x_r)$ itself is feasible in~\eqref{eq:optimal-prescriptor} remains open. 

We also emphasize that the strong LDP portrayed in Theorem~\ref{thm:sldp-iid-strong} has no continuous counterpart, which implies that the finite sample guarantees of Theorems~\ref{thm:pred:finite_sample} and~\ref{thm:presc:finite_sample} cannot be generalized.

\newpage

\begin{APPENDICES}
  \section{Proofs}
  \label{sec:proofs}

\proof{Proof of Theorem \ref{thm:sldp-iid}.}
Let $i_t\in\Xi$ be a particular realization of the random variable $\xi_t$ for each $t=1,\ldots,T$, and denote by $\mb P'$ the realization of the estimator $\hat{\mb P}_T$ corresponding to the sequence $\{i_t\}_{t=1}^T$. The probability of observing this sequence (in the given order) under model $\mb P$ can be expressed in terms of $\mb P'$ as
  \begin{equation}
  \label{eq:path-prob}
    \mb P^\infty (\xi_1=i_1, \dots, \xi_T=i_T) = \prod_{i\in\Xi}\mb P(i)^{T\mb P'(i)} = e^{T \sum_{i\in \Xi} \mb P'(i) \log \mb P(i)}.
  \end{equation}
  Set $\mc P_T=\mc P\cap \{0, 1/T, \dots, (T-1)/T,1\}^d$ and note that $\mb P^\infty(\hat {\mb P}_T \in \mc P_T) = 1$. By construction, the cardinality of $\mc P_{T}$ is bounded above by $(T+1)^d$.
   
  In the following, we denote the set of all sample paths in $\Xi^T$ that give rise to the same empirical distribution $\mb P'\in \mc P_T$ by $C_T(\mb P')$. The cardinality of $C_T(\mb P')$ coincides with the number of sample paths that visit state $i$ exactly $T\cdot \mb P'(i)$ times for all $i\in\Xi$, that is, we have
  \[
    \abs{C_T(\mb P')} = \frac{T!}{\prod_{i\in\Xi}(T\cdot \mb P'(i))!}.
  \]
  Stirling's approximation for factorials allows us to bound the cardinality of $C_T(\mb P')$ from both sides in terms of the entropy $H(\mb P')=-\sum_{i=1}^d \mb P'(i) \log \mb P'(i)$ of the empirical distribution $\mb P'$, that is,
  \begin{equation}
  \label{eq:stirling}
    (T+1)^{-d} e^{T H(\mb P')}\leq \abs{C_T(\mb P')} \leq e^{T H(\mb P')}.
  \end{equation}
  An elementary proof of these inequalities that does not involve Stirling's approximation is given by~\citet[Theorem~12.1.3]{cover2006elements}.
  
  Select an arbitrary Borel set $\mc D\subseteq\mc P$. For any $T\in\mb N$, we thus have
  \begin{align}
    \mb P^\infty ( \hat {\mb P}_T \in \mc D )  = &~ \sum_{\mb P'\in \mc D \cap \mc P_T} \mb P^\infty(\hat {\mb P}_T = \mb P')\nonumber \\
    \leq & ~(T+1)^d \cdot\max_{\mb P'\in \mc D\cap \mc P_T}  \mb P^\infty(\hat {\mb P}_T = \mb P') \nonumber \\
    \leq & ~(T+1)^d \cdot\max_{\mb P'\in \mc D\cap \mc P_T} \abs{C_T(\mb P')} e^{T \sum_{i\in \Xi} \mb P'(i) \log \mb P(i)} \nonumber \\
    \leq & ~(T+1)^d \cdot e^{-T \min_{\mb P'\in \mc D\cap \mc P_T} \D{\mb P'}{\mb P}} \nonumber \\
    \leq & ~(T+1)^d \cdot e^{-T \inf_{\mb P'\in \mc D} \D{\mb P'}{\mb P}}, \nonumber
  \end{align}
  where the first inequality exploits the estimate $\abs{\mc P_T} \leq (T+1)^d$, the second inequality holds due to~\eqref{eq:path-prob} and the definition of $C_T(\mb P')$, and the third inequality follows from the upper estimate in~\eqref{eq:stirling}. Taking logarithms on both sides of the above expression and dividing by $T$ yields
   \begin{align}
   \label{eq:sanov:ub}
    \frac 1 T \log \mb P^\infty ( \hat {\mb P}_T \in \mc D ) \leq & ~\frac{d\log (T+1)}{T} - \inf_{\mb P'\in \mc D} \D{\mb P'}{\mb P}.
  \end{align}
  Note that the finite sample bound~\eqref{eq:sanov:ub} does not rely on any properties of the set $\mc D$ besides measurability. The asymptotic upper bound~\eqref{eq:ldp_exponential_rates_ub} is obtained by taking the limit superior as $T$ tends to infinity on both sides of~\eqref{eq:sanov:ub}.

As for the lower bound~\eqref{eq:ldp_exponential_rates_lb}, recall that $\D{\mb P'}{\mb P}$ is continuous in $\mb P'$ as $\mb P>0$, see Proposition~\ref{prop:relative-entropy}(iii), and note that $\bigcup_{T\in\mb N} \mc P_T$ is dense in $\interior\mc D$. Thus, there exists $T_0\in\mb N$ and a sequence of distributions $\mb P'_T\in \mc P_T$, $T\in\mb N$, such that $\mb P'_T\in\interior \mc D$ for all $T \geq T_0$ and
\begin{equation}
	\label{eq:inf-int}
	\inf_{\mb P' \in \interior \mc D} \, \D{\mb P'}{\mb P} = \liminf_{T\rightarrow\infty} \, \D{\mb P'_T}{\mb P}.
\end{equation}
Fix any $T\geq T_0$ and let $(i_1,\ldots, i_T)$ be a sequence of observations that generates $\mb P'_T$. Then, we have
  \begin{align*}
    \label{eq:sanov:ub}
    \mb P^\infty ( \hat {\mb P}_T \in \mc D ) \geq &~  \mb P^\infty ( \hat {\mb P}_T =\mb P'_T) \\ 
    = &~ \abs{C_T(\mb P'_T)} \cdot \mb P^\infty (\xi_1=i_1, \dots, \xi_T=i_T)\\
    \geq &~ (T+1)^{-d} \cdot e^{TH(\mb P'_T)} \cdot e^{T\sum_{i\in\Xi} \mb P'_T(i)\log \mb P(i)}\\
    = &~ (T+1)^{-d} e^{-T \D{\mb P'_T}{\mb P}},
  \end{align*}
  where the first inequality holds because $\mb P_T'\in\interior\mc D\subseteq \mc D$, while the second inequality follows from~\eqref{eq:path-prob} and the lower estimate in~\eqref{eq:stirling}. This implies that
  \begin{equation*}
    \frac{1}{T} \log \mb P^\infty ( \hat {\mb P}_T \in \mc D ) \geq -\frac{d \log(T+1)}{T} - \D{\mb P'_T}{\mb P} \quad \forall T\geq T_0.
  \end{equation*}
 Taking the limit inferior as $T$ tends to infinity on both sides of the above inequality and using~\eqref{eq:inf-int} yields the postulated lower bound~\eqref{eq:ldp_exponential_rates_lb}. This completes the proof.  \hfill $\square$
\endproof

  \proof{Proof of Proposition~\ref{thm:strong-dual-final-finite}.}
  Applying \citep[Corollary~1]{ben2013robust} to the Burg entropy yields
  \begin{equation}
    \label{dual-finite-two}
    \hat c_{r}(x,\mb P') = \inf_{\alpha\geq \bar\gamma(x),\,\nu\geq 0} ~\sum_{i\in \Xi} \nu \log\left(\frac{\nu}{\alpha-\gamma(x, i)}\right) \mb P'(i) + \alpha + \nu (r-1)
  \end{equation}
  for any $r>0$. 
In the following we denote the objective function of~\eqref{dual-finite-two} by $g(\alpha, \nu)$, which is lower semicontinuous due to our conventions for the logarithm. As $g(\alpha,0)=\alpha$ and $\lim_{\nu\to\infty}g(\alpha, \nu)=\infty$, there must exist $\nu^\star(\alpha)\in \arg\min_{\nu\geq 0} g(\alpha, \nu)$ for any $\alpha\ge \bar \gamma(x)$. Indeed, if there is $i\in\Xi$ with $\alpha=\gamma(x,i)$ and $\mb P'(i)>0$, then $\nu^\star(\alpha)=0$. Otherwise, $\nu^\star(\alpha)$ is the unique solution of the first-order optimality condition
  \begin{align*}
    & \displaystyle\sum_{i\in\Xi}\left(\log\left( \frac{\nu^\star(\alpha)}{\alpha-\gamma(x, i)}\right) + 1\right) \mb P'(i) + r-1=0\\[0.5em]
    \iff ~& \displaystyle\sum_{i\in\Xi} \log\left( \frac{\nu^\star(\alpha)}{\alpha-\gamma(x, i)}\right) \mb P'(i) +r=0 \\
    \iff ~& \displaystyle \nu^\star(\alpha) = \exp\left(\sum_{i\in\Xi} \log\left(\alpha-\gamma(x, i)\right) \mb P'(i)- r \right).
  \end{align*}
  Thus, the partial minimum 
  \begin{equation}
  \label{eq:partial-minimum}
	g(\alpha, \nu^\star(\alpha))= \alpha - \exp\left( \sum_{i\in\Xi} \log\left(\alpha-\gamma(x, i)\right) \mb P'(i)-r \right)
	=\alpha - e^{-r} \prod_{i\in \Xi} \left(\alpha-\gamma(x, i) \right)^{\mb P'(i)}
   \end{equation}
   is readily recognized as the objective function of problem~\eqref{eq:strong-dual-final-finite}, which is manifestly continuous in $\alpha$ and inherits convexity from $g(\alpha,\nu)$. By using Jensen's inequality to interchange the logarithm and the expectation with respect to $\mb P'$ in the second expression in \eqref{eq:partial-minimum}, we find
  \begin{align}
  	\label{eq:lower-bound-g}
    g(\alpha, \nu^\star(\alpha)) \geq \alpha - \sum_{i\in \Xi}\left(\alpha-\gamma(x, i)\right) \mb P'(i) e^{-r}  \geq \alpha (1-e^{-r})+ e^{-r}  \sum_{i\in \Xi} \gamma(x, i)\, \mb P'(i).
  \end{align}
  As $r>0$, the above estimate implies that the partial minimum $g(\alpha, \nu^\star(\alpha))$ tends to infinity as $\alpha$ grows. Recalling that $\alpha$ is required to exceed $\bar\gamma(x)$, we may thus conclude that there exists $\alpha^\star$, such that $(\alpha^\star, \nu^\star(\alpha^\star))$ attains the minimum in~\eqref{eq:strong-dual-final-finite}. Moreover, as $\alpha\ge \bar\gamma(x)\ge \gamma(x,i)$ for all $i\in\Xi$, we have $g(\alpha^\star, \nu^\star(\alpha^\star))\leq \bar \gamma(x)$. Combining this upper bound with the lower bound~\eqref{eq:lower-bound-g} yields the estimate
  \begin{align*}
    \alpha^\star  \leq \frac{\bar \gamma(x)}{1-e^{-r}} -  \frac{e^{-r}}{1-e^{-r}} \sum_{i\in\Xi} \gamma(x, i) \,\mb P'(i)= \frac{\bar\gamma(x)-e^{-r} c(x, \mb P')}{1-e^{-r}}.
  \end{align*}
  Thus, the claim follows. %which can be added as a constraint without affecting the optimal value of~\eqref{eq:strong-dual-final-finite}.
  \hfill $\square$
  \endproof

\proof{Proof of Lemma \ref{lemma:r=0}.}
  As the cost function $\gamma(x,\xi)$ is continuous and its domain $X\times\Xi$ is compact, the Heine-Cantor theorem implies that $\gamma(x,\xi)$ is uniformly continuous on its domain. Consider now an arbitrary converging sequence $(x_i,\mb P_i)$, $i\in\mb N$, in $X\times\mc P$, and denote its limit by $(x, \mb P)$. The uniform continuity of the cost function ensures that for every $\delta>0$ there exists $N_\delta\in\mb N$ such that $\abs{\gamma(x_i, \xi) -\gamma(x,\xi)}\leq \delta$ uniformly across all $\xi\in\Xi$ and $i\ge N_\delta$, which in turn implies that 
\[
  \abs{\int_{\Xi} \gamma(x, \xi) \mb P( \d\xi)- \int_{\Xi} \gamma(x_i, \xi) \d \mb P_i(\xi)}  \leq  \abs{\int_{\Xi}\gamma(x, \xi)\d \mb P(\xi) - \int_{\Xi} \gamma(x, \xi) \d \mb P_i(\xi)} + \delta
  \quad \forall i\ge N_\delta.
\]
As the integrant $\gamma( x, \xi)$ on the right hand side of the above inequality is continuous and bounded in~$\xi$, and as the sequence $\mb P_i$, $i\in\mb N$, converges weakly to $\mb P$, we thus have
\(
  \lim_{i\to \infty}\abs{\int_{\Xi} \gamma(x, \xi) \mb P(\d\xi) - \int_{\Xi} \gamma(x_i, \xi)  \mb P_i(\d\xi)} \leq \delta.
\)
As $\delta>0$ was chosen arbitrary, we may finally conclude that
\[
  \lim_{i\to\infty}\int_{\Xi} \gamma(x_i, \xi) \mb P_i (\d\xi)= \int_{\Xi} \gamma(x, \xi) \mb P(\d\xi).
\]
The claim then follows because the converging sequence $(x_i,\mb P_i)$, $i\in\mb N$, was chosen arbitrary. \hfill $\square$
\endproof

In order to prove Proposition~\ref{thm:strong-dual-final} we need three auxiliary lemmas. As a starting point, we first exploit the definition of the relative entropy between arbitrary distributions on a compact state space $\Xi$ to re-express the distributionally robust predictor explicitly as
\begin{equation}
  \label{eq:primal-problem}
  \begin{array}{rcl}
    \hat c_r(x, \mb P')=& \displaystyle \sup_{\mb P\in \mc P} & \displaystyle\int_{\Xi} \gamma(x, \xi)\, \d\mb P(\xi)\\[1em]
    &\st & \mb P' \ll \mb P \\[0.5em]
                            && \displaystyle\int_{\Xi} \log \left( \frac{\d \mb P'}{\d\mb P} (\xi)\right) \d \mb P'(\xi) \leq r.
  \end{array}
\end{equation}
Under the standing assumptions that $X$ and $\Xi$ are compact, while the cost function $\gamma(x, \xi)$ is jointly continuous in its arguments, the feasible set of problem \eqref{eq:primal-problem} can be restricted to distributions $\mb P$ that are absolutely continuous with respect to $\mb P'$ except perhaps on the compact set $\Xi^\star(x)=\arg\max_{\xi\in\Xi}\gamma(x,\xi)$. 

\begin{lemma}[Absolutely continuous representation of $\hat c_r$]
  \label{thm:primal:continuous}
  If $r\ge 0$ and $\bar \gamma(x)=\max_{\xi\in\Xi}\gamma (x,\xi)$ denotes the worst-case cost function, then the distributionally robust predictor $\hat c_r$ admits the equivalent representation
  \begin{equation}
    \label{eq:primal-problem:representation}
    \begin{array}{rcl}
      \hat c_{r}(x, \mb P')=& \displaystyle \sup_{\substack{\mb P_c\in\mc P\\ p\in[0,1]}} & \displaystyle p\cdot \int_{\Xi} \gamma(x, \xi)\, \d\mb P_c(\xi) + (1-p)\cdot \bar\gamma(x) \\
      & \st & \mb P'\ll p\cdot \mb P_c \ll \mb P'\\[0.5em]
      && \displaystyle \int_{\Xi} \log \left( \frac{1}{p}\cdot \frac{\d \mb P'}{\d\mb P_c} (\xi)\right) \d\mb P'(\xi) \leq r.
    \end{array}
  \end{equation}
\end{lemma}

\proof{Proof.}
We first show that \eqref{eq:primal-problem} provides an upper bound on \eqref{eq:primal-problem:representation}. To this end, choose any $\mb P_c$ and $p$ feasible in~\eqref{eq:primal-problem:representation}, and define $\mb P =p\cdot \mb P_c + (1-p)\cdot \delta_{\xi^\star}\in \mc P$, where $\xi^\star\in \Xi^\star(x)$ represents an arbitrary worst-case scenario.

Note that $p>0$ for otherwise $ \mb P'\not \ll p\cdot \mb P_c$. The constraints of~\eqref{eq:primal-problem:representation} and the construction of $\mb P$ thus imply that $\mb P'\ll \mb P_c \ll \mb P$. By the Radon-Nikodym theorem \citep{nikodym_sur_1930}, we then have 
\begin{align*}
	\mb P'[A] & = \int_A \frac{\d \mb P'}{\d\mb P}(\xi)\, \d\mb P(\xi)
\end{align*}
and
\begin{align*}
	\mb P'[A] & = \int_A \frac{\d \mb P'}{\d\mb P_c}(\xi)\, \d\mb P_c(\xi) = \int_A\frac{1}{p}\cdot \frac{\d \mb P'}{\d\mb P_c}(\xi) \cdot p\, \d\mb P_c(\xi)\\
	&= \int_A\frac{1}{p}\cdot \frac{\d \mb P'}{\d\mb P_c}(\xi) \, \d\mb P(\xi) - \frac{1-p}{p}\cdot  \frac{\d \mb P'}{\d\mb P_c}(\xi^\star) \le \int_A\frac{1}{p}\cdot \frac{\d \mb P'}{\d\mb P_c}(\xi) \, \d\mb P(\xi)
\end{align*}
for all Borel sets $A\subseteq \Xi$, where the inequality in the last expression holds because Radon-Nikodym derivatives are pointwise non-negative. In summary, the above reasoning implies that
\[
	 \frac{\d \mb P'}{\d\mb P}(\xi) \le \frac{1}{p}\cdot \frac{\d \mb P'}{\d\mb P_c}(\xi)  \quad \mb P\text{-a.s.}
\]
In fact, as $\mb P'\ll\mb P$, the above inequality even holds almost surely with respect to $\mb P'$, which in turn implies
  \begin{align*}
    \int_{\Xi} \log \left( \frac{\d \mb P'}{\d\mb P} (\xi)\right) \d\mb P'(\xi)  
    \leq \int_{\Xi} \log \left( \frac{1}{p}\cdot \frac{\d \mb P'}{\d\mb P_c} (\xi)\right) \d\mb P'(\xi)\leq r.
  \end{align*}
Thus, $\mb P$ is feasible in~\eqref{eq:primal-problem}. Moreover, it is easy to verify that the objective value of $\mb P$ in~\eqref{eq:primal-problem} is equal to that of $(\mb P_c, p)$ in \eqref{eq:primal-problem:representation}. Thus, \eqref{eq:primal-problem} provides an upper bound on \eqref{eq:primal-problem:representation}.

  It remains to be shown that~\eqref{eq:primal-problem} provides also a lower bound on~\eqref{eq:primal-problem:representation}. To this end, choose any $\mb P$ feasible in~\eqref{eq:primal-problem}, define $\Xi^{+}=\{\xi\in\Xi:\d\mb P'/\d\mb P(\xi)>0\}$ as the (measurable) event in which the Radon-Nikodym derivative of $\mb P'$ with respect to $\mb P$ is strictly positive, and set $p=\mb P[\Xi^{+}]$. Note that $p>0$, for otherwise the relations $\mb P[\Xi^+]=0$ and $\mb P'\ll\mb P$ would imply via the Radon-Nikodym theorem that 
  \[
  0= \mb P'[\Xi^+]=\int_{\Xi^+} \frac{\d\mb P'}{\d\mb P} (\xi)\,\d\mb P(\xi)=\int_{\Xi} \frac{\d\mb P'}{\d\mb P} (\xi)\,\d\mb P(\xi) =\mb P'[\Xi]=1.
  \] 
  Next, we define $\mb P_c\in\mc P$ through $\mb P_c[A]=\mb P[A\cap \Xi^+]/p$ for all Borel sets $A\subseteq\Xi$. By construction, $\mb P'$ is absolutely continuous with respect to $\mb P_c$. To see this, note that for any Borel set $A\subseteq \Xi$ we have
  \begin{align*}
  	\mb P_c[A] = 0 ~ \iff ~ \mb P[A\cap \Xi^+]=0 ~ \implies ~ \mb P'[A\cap \Xi^+]=0 ~\iff \mb P'[A]=0,
  \end{align*}
  where the implication holds because $\mb P'\ll\mb P$. Conversely, one can also show that $\mb P_c$ is absolutely continuous with respect to $\mb P'$. To see this, assume that $\mb P_c[A]>0$ for some Borel set $A\subseteq \Xi$. By the construction of $\mb P_c$ we thus have $\mb P[A\cap \Xi^+]>0$. The Radon-Nikodym theorem further implies that
  \[
  	\mb P'[A]= \int_A \frac{\d\mb P'}{\d\mb P} (\xi)\,\d\mb P(\xi)= \int_{A\cap \Xi^+} \frac{\d\mb P'}{\d\mb P} (\xi)\,\d\mb P(\xi)>0,
  \]
  where the inequality holds because the integral of a strictly positive function over a set of strictly positive measure must be strictly positive. We have thus shown that $\mb P_c[A]>0$ implies $\mb P'[A]>0$ or, by contraposition, that $\mb P'[A]=0$ implies $\mb P_c[A]=0$. In summary, the above reasoning ensures that $\mb P'\ll p\cdot \mb P_c \ll \mb P'$.
  
Using the Radon-Nikodym theorem once again, we have for all Borel sets $A\subseteq \Xi$ that
\begin{align*}
	\int_A \frac{\d \mb P'}{\d\mb P_c}(\xi) \, \d\mb P_c(\xi) = \mb P'[A] =\int_A \frac{\d \mb P'}{\d\mb P}(\xi)\, \d\mb P(\xi)=  
	\int_A \frac{\d \mb P'}{\d\mb P}(\xi)\cdot p \, \d\mb P_c(\xi) ,
\end{align*}
where the last equality holds because $\d\mb P'/\d\mb P(\xi)=0$ for all $\xi\notin\Xi^+$ and because $\mb P= p\cdot \mb P_c$ when restricted to the Borel $\sigma$-algebra on $\Xi^+$. As the above equality holds for all  Borel sets $A\subseteq \Xi$, we find
\[
	\frac{\d \mb P'}{\d\mb P_c}(\xi)=p\cdot \frac{\d \mb P'}{\d\mb P}(\xi)\quad \mb P_c\text{-a.s.} \qquad \implies\qquad 
	\frac{1}{p}\cdot \frac{\d \mb P'}{\d\mb P_c}(\xi) = \frac{\d \mb P'}{\d\mb P}(\xi)\quad \mb P'\text{-a.s.},
\]
where the implication holds because $p>0$ and $\mb P'\ll\mb P_c$. The last identity and our standard convention $0\log(0)=0$ ensure that
\[
	\int_{\Xi} \log \left( \frac{1}{p}\cdot \frac{\d \mb P'}{\d\mb P_c} (\xi)\right) \d\mb P'(\xi) = \int_{\Xi} \log \left( \frac{\d \mb P'}{\d\mb P} (\xi)\right) \d \mb P'(\xi) \leq r.
\]
 Thus, $(p, \mb P_c)$ is feasible in~\eqref{eq:primal-problem:representation}. 
 
 Assume now that $p<1$, and define $\mb P_\perp\in\mc P$ through $\mb P_\perp[A]=\mb P[A\backslash \Xi^+]/(1-p)$ for all Borel sets $A\subseteq \Xi$. By construction, $\mb P_\perp$ is thus singular with respect to $\mb P_c$, and we have $\mb P=p\cdot \mb P_c+(1-p) \cdot \mb P_\perp$. This implies that
 \[
 	\int_\Xi \gamma (x,\xi)\, \d \mb P(\xi) = p\cdot \int_{\Xi} \gamma(x, \xi)\, \d\mb P_c(\xi) + (1-p)\cdot \int_\Xi \gamma(x,\xi)\,\d\mb P_\perp(\xi) 
	\le p\cdot \int_{\Xi} \gamma(x, \xi)\, \d\mb P_c(\xi) + (1-p)\cdot \bar\gamma(x).
 \]
 If $p=1$, on the other hand, we trivially have $\mb P=\mb P_c$ and
\[
	\int_\Xi \gamma (x,\xi)\, \d \mb P(\xi) 	= p\cdot \int_{\Xi} \gamma(x, \xi)\, \d\mb P_c(\xi) + (1-p)\cdot \bar\gamma(x).
\]
 In summary, we have thus shown that for every $\mb P$ feasible in~\eqref{eq:primal-problem} there exists $(p,\mb P_c)$ feasible in~\eqref{eq:primal-problem:representation} with the same or with a larger objective value. This implies that~\eqref{eq:primal-problem} provides a lower bound on~\eqref{eq:primal-problem:representation}. \hfill $\square$
  \endproof

  Define now the family of predictors
   \begin{equation}
    \label{eq:primal-problem:approximation}
    \begin{array}{rcl}
      \hat c_{r,\epsilon}(x, \mb P')=& \displaystyle \sup_{\substack{\mb P_c\in\mc P\\ p\in[0,1]}} & \displaystyle p\cdot \int_{\Xi} \gamma(x, \xi)\, \d\mb P_c(\xi) + (1-p)\cdot (\bar\gamma(x) +\epsilon)\\
      & \st & \mb P'\ll p\cdot \mb P_c \ll \mb P'\\[0.5em]
      && \displaystyle \int_{\Xi} \log \left( \frac{1}{p}\cdot \frac{\d \mb P'}{\d\mb P_c} (\xi)\right) \d\mb P'(\xi) \leq r
    \end{array}
  \end{equation}
 parameterized by $r> 0$ and $\epsilon\ge 0$. 
 We first show that $\hat c_{r, \epsilon}$ uniformly approximates $\hat c_{r}$. 

\begin{lemma}[Uniform approximation of $\hat c_r$]
  \label{lemm:uniform_approximation}
  If $r>0$ and $\epsilon\ge 0$, then 
  \[
    \hat c_r(x, \mb P') \leq \hat c_{r,\epsilon}(x, \mb P') \leq \hat c_r(x, \mb P') + \epsilon\quad \forall x\in X,~\mb P'\in \mc P.
   \]
\end{lemma}
\proof{Proof.}
The claim follows immediately by comparing the absolutely continuous representations for $\hat c_r$ derived in Lemma~\ref{thm:primal:continuous} with the definition of~$\hat c_{r,\epsilon}$ in~\eqref{eq:primal-problem:approximation}. \hfill $\square$
\endproof

Next, we demonstrate that the predictors~$\hat c_{r,\epsilon}$ defined in~\eqref{eq:primal-problem:approximation} admit a dual representation in the form of a univariate convex optimization problem.

\begin{lemma}[Dual representation of $\hat c_{r,\epsilon}$]
  \label{lemma:geom-mean}
  If $r>0$, $\epsilon\ge0$ and $\bar\gamma(x)=\max_{\xi\in \Xi} \, \gamma(x, \xi)$ denotes the worst-case cost function, then the predictor $\hat c_{r,\epsilon}$ defined in~\eqref{eq:primal-problem:approximation} satisfies
  \begin{equation}
    \label{eq:strong-dual-epsilon}
      \hat c_{r,\epsilon}(x,\mb P')\le  \min_{\alpha\ge\bar\gamma(x)+\epsilon} \alpha - e^{-r}\cdot \exp\left( \int_{\Xi} \log\left(\alpha-\gamma(x, \xi)\right) \d\mb P'(\xi)\right),
  \end{equation}
  and the problem on the right hand side of~\eqref{eq:strong-dual-epsilon} has a minimizer $\alpha^\star \leq \frac{\bar\gamma(x)+\epsilon -e^{-r} c(x, \mb P')}{1-e^{-r}}$. Moreover, if $\epsilon>0$, then~\eqref{eq:strong-dual-epsilon} becomes an equality.
 \end{lemma}

\proof{Proof.}
By applying the variable transformations $p\leftarrow 1-p$ and $\mb P_c\leftarrow p\cdot\mb P_c$, the non-convex optimization problem~\eqref{eq:primal-problem:approximation} can be reformulated as 
\begin{equation}
	\label{eq:primal-problem:approximation2}
    \begin{array}{rcl}
      \hat c_{r,\epsilon}(x,\mb P')= & \displaystyle \sup_{\substack{\mb P_c\ge 0\\ p\ge 0}} & \displaystyle \int_{\Xi} \gamma(x, \xi)\, \d\mb P_c(\xi) + p\cdot (\bar\gamma(x) +\epsilon)\\
      &\st & \mb P'\ll \mb P_c \ll \mb P'\\[0.3em]
      && \displaystyle \int_{\Xi} \log \left(\frac{\d \mb P'}{\d\mb P_c} (\xi)\right) \d\mb P'(\xi) \leq r\\[1.0em]
      && \displaystyle \int_{\Xi} \d\mb P_c(\xi) +p =1,
    \end{array}
\end{equation}
which is manifestly convex. The condition $\mb P_c\ge 0$ abbreviates the requirement that $\mb P_c$ is a finite non-negative Borel measure supported on $\Xi$. Note also that the normalization of $\mb P_c$ is now enforced through an explicit constraint. Next, we eliminate the decision variable $\mb P_c$ from~\eqref{eq:primal-problem:approximation2} by re-expressing it in terms of $\mb P'$ and the Radon-Nikodym derivative $\Lambda(\xi)=\d \mb P_c/\d\mb P'(\xi)$. In particular, as $\mb P_c\ll \mb P'$ and $\mb P'\ll \mb P_c$, we have $\d\mb P'/\d\mb P_c(\xi)=( \d \mb P_c/\d\mb P'(\xi))^{-1}=\Lambda(\xi)^{-1}$ almost everywhere with respect to $\mb P'$. 
Thus, problem~\eqref{eq:primal-problem:approximation2} is equivalent to
\begin{equation}
	\label{eq:primal-problem:approximation3}
    \begin{array}{rcl}
      \hat c_{r,\epsilon}(x,\mb P')= & \displaystyle \sup_{\substack{\Lambda\ge 0\\ p\ge 0}} & \displaystyle \int_{\Xi} \gamma(x, \xi) \Lambda(\xi)\, \d\mb P'(\xi) + p\cdot (\bar\gamma(x) +\epsilon)\\
      &\st & \displaystyle - \int_{\Xi} \log \left(\Lambda (\xi)\right) \d\mb P'(\xi) \leq r\\[2ex]
      && \displaystyle \phantom{-} \int_{\Xi} \Lambda(\xi)\, \d\mb P'(\xi) +p =1,
    \end{array}
\end{equation}
where the condition $\Lambda\ge 0$ abbreviates the requirement that $\Lambda$ is a non-negative Borel-measurable function on $\Xi$. The first (relative entropy) constraint ensures that $\Lambda(\xi)>0$ almost surely with respect to $\mb P'$, while the second (normalization) constraint ensures that the measure induced by $\Lambda$ and $\mb P'$ is finite. Thus, the constraint that $\mb P'$ be absolutely continuous with respect to $\mb P_c$, which is explicitly imposed in~\eqref{eq:primal-problem:approximation2}, remains implicitly enforced in~\eqref{eq:primal-problem:approximation3}. The Lagrangian dual of the convex maximization problem~\eqref{eq:primal-problem:approximation3} is given by
\begin{equation}
	\label{eq:dual-predictor-continuous}
	\inf_{\alpha\in\Re,\,\nu\in\Re_+} g(\alpha,\nu),
\end{equation}
where the dual objective function can be represented as
\[
	g(\alpha,\nu)=   \sup_{\substack{\Lambda\ge 0\\ p\ge 0}} \int_{\Xi} \left[\gamma(x, \xi)-\alpha\right] \Lambda(\xi) + \nu \log \left( \Lambda (\xi)\right) \, \d\mb P'(\xi) 
	+ \left(\bar\gamma(x)+\epsilon-\alpha\right) p + \alpha + \nu r.
\]
By weak duality, the dual problem~\eqref{eq:dual-predictor-continuous} provides an upper bound on~$\hat c_{r,\epsilon}(x,\mb P')$. Note that the supremum over $p\ge 0$ in the definition of $g(\alpha, \nu)$ is unbounded if $\bar\gamma(x)+\epsilon> \alpha$ and evaluates to 0 otherwise. Thus, the dual problem~\eqref{eq:dual-predictor-continuous} includes the implicit constraint $\bar\gamma(x)+\epsilon\le \alpha$. We henceforth assume that this constraint holds, and we assume that the logarithm of any nonpositive number is defined as~$-\infty$. Under this premise we may remove the redundant constraint $\Lambda\ge 0$ and reformulate the dual objective function as
\[
  \begin{array}{rl}
    g(\alpha, \nu) & =  \displaystyle \sup_{\Lambda} \int_{\Xi} \left[\gamma(x, \xi)- \alpha\right] \Lambda(\xi) + \nu \log\left( \Lambda (\xi)\right) \d\mb P'(\xi)  + \alpha + \nu r\\[2ex]
                   & = \displaystyle \int_{\Xi} \sup_{\lambda} \left[\gamma(x, \xi)-\alpha\right] \lambda + \nu \log\left( \lambda \right) \d\mb P'(\xi) + \alpha + \nu r \\[2ex]                  
                   & = \displaystyle \int_{\Xi} \nu \log\left(\frac{\nu}{\alpha-\gamma(x, \xi)}\right) \d\mb P'(\xi) + \alpha + \nu (r-1),
  \end{array}
\]
where the supremum over all measurable functions $\Lambda$ can be moved inside the integral and converted to a supremum over all scalars $\lambda$ by appealing to~\citep[Theorem~14.60]{rockafellar1998variational}. The third equality in the last line follows from an explicit solution of the convex maximization problem over~$\lambda$, which has a unique well-defined solution because $\gamma(x,\xi)-\alpha\ge \epsilon$ for all $\xi\in\Xi$. The dual problem~\eqref{eq:dual-predictor-continuous} is thus equivalent to
    \begin{equation}
    \label{eq:dual}
    \min_{\alpha\ge\bar\gamma (x)+\epsilon, \nu\ge 0} \int_{\Xi} \nu \log\left(\frac{\nu}{\alpha-\gamma(x, \xi)}\right) \d\mb P'(\xi) + \alpha + \nu (r-1),
    \end{equation}
  which constitutes a finite-dimensional convex minimization problem. If $\epsilon=0$, then~\eqref{eq:dual} can be viewed as a continuous version of~\eqref{dual-finite-two}. The dual objective function $g(\alpha,\nu)$ is lower semicontinuous. To see this, note that the function inside the integral in~\eqref{eq:dual} is lower semicontinuous in $(\alpha, \nu)$ due to our conventions for the logarithm and because lower semicontinuity is preserved under integration thanks to Fatou's lemma. Following a similar reasoning as in the proof of Proposition~\ref{thm:strong-dual-final-finite}, one can now show that the minimum of~\eqref{eq:dual} is always attained by some $\alpha^\star$ and $\nu^\star$ and that~\eqref{eq:dual} is equivalent to the one-dimensional convex program
  \begin{equation}
  	\label{eq:strong-dual-epsilon2}
      \min_{\alpha\ge\bar\gamma(x)+\epsilon} \alpha - e^{-r}\cdot \exp\left( \int_{\Xi} \log\left(\alpha-\gamma(x, \xi)\right) \d\mb P'(\xi)\right),
  \end{equation}
  which has a minimizer $\alpha^\star \leq \frac{\bar\gamma(x)+\epsilon -e^{-r} c(x, \mb P')}{1-e^{-r}}$. Details are omitted for the sake of brevity.  
  
  Note that~\eqref{eq:strong-dual-epsilon2} coincides with coincides with the optimization problem on the right hand side of~\eqref{eq:strong-dual-epsilon}. The above arguments thus imply that~\eqref{eq:strong-dual-epsilon2} provides an upper bound on $\hat c_{r,\epsilon}(x,\mb P')$. To prove that this upper bound is in fact exact for $\epsilon>0$, it remains to be shown that the duality gap between~\eqref{eq:primal-problem:approximation3} and~\eqref{eq:dual} vanishes. We will do so by constructing a pair of primal and dual feasible solutions whose objective values coincide.

For $\epsilon>0$ the minimization problem~\eqref{eq:dual} satisfies Slater's constraint qualification because its convex objective function is finite-valued on the feasible set and because the decision variables are only subject to lower bounds. Thus, $(\alpha^\star, \nu^\star)$ solves \eqref{eq:dual} if and only if it satisfies the Karush-Kuhn-Tucker (KKT) conditions
  \begin{equation*}
    %\label{eq:kkt}
    \begin{array}{l@{\qquad}l}
        \alpha \geq \bar\gamma (x) + \epsilon, ~ \nu \geq 0 & \text{(primal feasibility)} \\[0.5em]
        \lambda\geq 0 & \text{(dual feasibility)} \\[0.5em]
        \int_{\Xi} \frac{\nu}{\alpha-\gamma(x, \xi)}\d \mb P'(\xi)+\lambda = 1 & \text{(stationarity with respect to~$\alpha$)} \\[0.5em]
        \int_{\Xi} \log\left(\frac{\alpha-\gamma(x, \xi)}{\nu}\right)\d \mb P'(\xi) = r & \text{(stationarity with respect to~$\nu$)} \\[0.5em]
        \lambda\cdot (\alpha - \bar\gamma (x)-\epsilon) = 0, & \text{(complementary slackness)} 
      \end{array}
  \end{equation*}
  where $\lambda$ represents the Lagrange multiplier of the constraint $\alpha\geq \bar\gamma (x) + \epsilon$. Given any solution $(\alpha^\star, \nu^\star,\lambda^\star)$ of the KKT conditions, we can now introduce a Borel-measurable function $\Lambda^\star(\xi)= \frac{\nu^\star}{\alpha^\star-\gamma(x, \xi)}$ and define $p^\star = \lambda^\star$. As $\alpha^\star\ge \bar\gamma(x)+\epsilon$, the function $\Lambda^\star$ is strictly positive on $\Xi$. The two stationarity conditions thus imply that $(\Lambda^\star,p^\star)$ is feasible in~\eqref{eq:primal-problem:approximation3}. Moreover, we have
  \begin{align*}
  	g(\alpha^\star, \nu^\star)& = \int_{\Xi} \nu^\star \log\left(\frac{\nu^\star}{\alpha^\star-\gamma(x, \xi)}\right) \d\mb P'(\xi) + \alpha^\star + \nu^\star (r-1)  = \alpha^\star-\nu^\star \\
	& = \alpha^\star+ \int_\Xi (\gamma(x,\xi)-\alpha^\star) \frac{\nu^\star}{\alpha^\star-\gamma(x,\xi)} \, \d\mb P'(\xi) \\
	& = \alpha^\star+ \int_\Xi \gamma(x,\xi) \Lambda^\star(\xi) \, \d\mb P'(\xi) -\alpha^\star\cdot(1-\lambda^\star)\\
	& =  \int_{\Xi}\gamma(x, \xi) \Lambda^\star(\xi) \, \d \mb P'(\xi) + p^\star\cdot \left(\bar\gamma(x)+ \epsilon\right)
  \end{align*}
where the first equality follows from the definition of $g$, the second equality holds due to the stationarity condition for $\nu$, the fourth equality follows from the definition of $\Lambda^\star$ and the stationarity condition for $\alpha$, and the last equality exploits the complementary slackness condition as well as the definition of~$p^\star$. 

We have thus shown that the objective value of~$(\Lambda^\star,p^\star)$ in~\eqref{eq:primal-problem:approximation3} coincides with the objective value of~$(\alpha^\star, \nu^\star)$ in~\eqref{eq:dual}, which certifies that the duality gap between~\eqref{eq:primal-problem:approximation3} and~\eqref{eq:dual} vanishes. \hfill $\square$
\endproof

Armed with Lemmas~\ref{lemm:uniform_approximation} and~\ref{lemma:geom-mean}, we are now ready to prove Proposition~\ref{thm:strong-dual-final}.

\proof{Proof of Proposition~\ref{thm:strong-dual-final}.}
For every $\epsilon>0$ we have
\begin{align*}
	\hat c_r(x,\mb P')& \displaystyle \le  \min_{\alpha\ge\bar\gamma(x)} \alpha - e^{-r}\cdot \exp\left( \int_{\Xi} \log\left(\alpha-\gamma(x, \xi)\right) \d\mb P'(\xi)\right) \\
	& \displaystyle \le  \min_{\alpha\ge\bar\gamma(x)+\epsilon} \alpha - e^{-r}\cdot \exp\left( \int_{\Xi} \log\left(\alpha-\gamma(x, \xi)\right) \d\mb P'(\xi)\right)  = \hat c_{r,\epsilon}(x,\mb P')\le \hat c_r(x,\mb P')+\epsilon,
\end{align*}
where the first inequality follows from Lemma~\ref{lemma:geom-mean} for $\epsilon=0$, the equality follows from Lemma~\ref{lemma:geom-mean} for $\epsilon>0$, and the last inequality follows from Lemma~\ref{lemm:uniform_approximation}. As the above inequalities remain valid for all $\epsilon >0$, we may conclude that~\eqref{eq:strong-dual} holds. The claim that the minimization problem on the right hand side of~\eqref{eq:strong-dual} has an optimizer $\alpha^\star\le\frac{\bar\gamma(x) -e^{-r} c(x, \mb P')}{1-e^{-r}}$ follows from Lemma~\ref{lemma:geom-mean}. \hfill $\square$
\endproof

\proof{Proof of Proposition~\ref{prop:pred:feasibility:revisited}.} Fix $\epsilon>0$ and recall from Lemma~\ref{lemma:geom-mean} that $\hat c_{r,\epsilon}(x,\mb P')$ coincides with the optimal value of the univariate convex minimization problem on the right hand side of~\eqref{eq:strong-dual-epsilon}. Throughout this proof we assume without loss of generality that this problem accommodates the extra constraint $\alpha \le\frac{\bar\gamma(x)+\epsilon -e^{-r} c(x, \mb P')}{1-e^{-r}}$. Indeed, Lemma~\ref{lemma:geom-mean} guarantees that this constraint has no impact on the problem's optimal value.

In the first part of the proof we demonstrate that the compactified feasible set
\[
	\left\{ \alpha\in\Re~ :~ \bar \gamma(x)+\epsilon\le \alpha \le \frac{\bar\gamma(x)+\epsilon -e^{-r} c(x, \mb P')}{1-e^{-r}} \right\}
\]
 of problem~\eqref{eq:strong-dual-epsilon} with the redundant upper bound on $\alpha$ represents a continuous set-valued mapping parameterized in $(x,\mb P')$. To this end, we note that the worst-case cost $\bar\gamma(x)$ is continuous in $x$ by Berge's maximum theorem \citep[pp.~115--116]{berge1963topological}, which applies because $\Xi$ is compact and $\gamma(x,\xi)$ is jointly continuous in $x$ and $\xi$. Lemma~\ref{lemma:r=0} further implies that~$c(x, \mb P')$ is continuous in $x$ and $\mb P'$. As both the upper and the lower bound on $\alpha$ depend continuously on $(x,\mb P')$, we conclude that the feasible set mapping is indeed continuous.

  In the second part of the proof we argue that the objective function of~\eqref{eq:strong-dual-epsilon} is continuous in $(\alpha, x,\mb P')$. To this end, recall that the cost function $\gamma(x,\xi)$ is uniformly continuous on its compact domain $X\times \Xi$. Consider now an arbitrary converging sequence $(\alpha_i, x_i,\mb P'_i)$, $i\in\mb N$, in $\mb R\times X\times\mc P$ such that $\alpha_i\ge \bar\gamma(x_i)+\epsilon$ for all $i\in\mb N$, and denote its limit by $(\alpha, x, \mb P')$. The uniform continuity of the cost function ensures that for every $\delta>0$ there exists $N_\delta\in\mb N$ such that $|\alpha_i-\alpha|\le \delta$ and $\abs{\gamma(x_i, \xi) -\gamma(x,\xi)}\leq \delta$ uniformly across all $\xi\in\Xi$ and $i\ge N_\delta$. As the natural logarithm is Lipschitz continuous on $[\epsilon,\infty)$ with Lipschitz constant $1/\epsilon$, we thus have
  \[
  	\abs{\log(\alpha_i-\gamma(x_i, \xi)) -\log( \alpha-\gamma( x,\xi))}\leq {2\delta}/{\epsilon}\quad\forall \xi\in\Xi,~i\in\mb N.
\]
This implies that
\begin{align*}
  & \abs{\int_{\Xi} \log(\alpha_i-\gamma(x_i, \xi)) \d \mb P'_i - \int_{\Xi} \log(\alpha-\gamma( x, \xi)) \d {\mb P}' } \\ & \qquad \leq  \abs{ \int_{\Xi} \log(\alpha-\gamma( x, \xi)) \d \mb P'_i - \int_{\Xi} \log(\alpha-\gamma( x, \xi)) \d{\mb P}' } + 2 \delta/\epsilon.
\end{align*}
As the limit of the chosen sequence satisfies $\alpha-\gamma(x, \xi)\geq \epsilon >0$, the integrand on the right hand side of the above inequality is continuous and bounded in $\xi$. By the definition of weak convergence we thus find
\[
  \lim_{i\to \infty}\abs{\int_{\Xi} \log(\alpha_i-\gamma(x_i, \xi)) \d \mb P'_i - \int_{\Xi} \log(\alpha-\gamma( x, \xi)) \d {\mb P}' } \leq 2 \delta/\epsilon.
\]
As $\delta>0$ was chosen arbitrary, it follows that $\int_{\Xi} \log(\alpha_i-\gamma(x_i, \xi)) \d \mb P_i$ converges to $\int_{\Xi} \log(\bar\alpha-\gamma(\bar x, \xi)) \d \bar{\mb P}'$, which establishes that the objective function of problem~\eqref{eq:strong-dual-epsilon} is continuous in $(\alpha, x, \mb P')$.

In summary, we have shown that the compactified feasible set of problem~\eqref{eq:strong-dual-epsilon} is continuous in $(x,\mb P')$ and that the objective function of~\eqref{eq:strong-dual-epsilon} is continuous in $(\alpha, x, \mb P')$. Thus, $\hat c_{r,\epsilon} (x, \mb P')$ is continuous by Berge's maximum theorem \citep{berge1963topological}. As $\hat c_{r,\epsilon}(x, \mb P')$ uniformly approximates $\hat c_r(x, \mb P')$ for $\epsilon\downarrow 0$ (see Lemma~\ref{lemm:uniform_approximation}), and as uniform limits of continuous functions are continuous, we conclude that $\hat c_r(x, \mb P')$ is continuous. \hfill $\square$
\endproof

\proof{Proof of Theorem \ref{thm:pred:feasibility-revisited}.}
We first establish feasibility for $r\ge 0$. From Proposition~\ref{prop:pred:feasibility:revisited} we already know that $\hat c_r$ is continuous on $X\times\mc P$, that is, $\hat c_r\in\mc C$. To show that the out-of-sample disappointment of $\hat c_r$ decays sufficiently fast, we fix any decision $x\in X$ and an arbitrary distribution $\mb P_0\in\mc P$, and we define 
  \[
    \mc D(x, \mb P_0)= \left\{ \mb P'\in\mc P : c(x, \mb P_0) >  \hat c_{r}(x, \mb P') \right\}\quad \text{and}\quad
    \bar{\mc D}(x, \mb P_0)= \left\{ \mb P'\in\mc P : c(x, \mb P_0) \ge  \hat c_{r}(x, \mb P') \right\}
  \]
as the corresponding strict and weak disappointment sets. Using this notation, we need to demonstrate that the probability of the event $\hat {\mb P}_T\in \mc D(x, \mb P_0)$ decays at a rate of at least $r$. We will prove this assertion by case distinction depending on the probability of the set of worst-case scenarios, $\Xi^\star(x)= \arg \max_{\xi\in \Xi} \gamma(x, \xi)$, which is non-empty and compact because $\gamma$ is continuous and $\Xi$ is compact.

 \paragraph{Case 1:} Assume first that $\mb P_0(\Xi^\star(x))=1$. As the support of $\hat{\mb P}_T$ is $\mb P^\infty$-almost surely a subset of the support of $\mb P_0$, we may thus conclude that $\hat {\mb P}_T$ is $\mb P_0^\infty$-almost surely supported on $\Xi^\star(x)$. Setting $\bar \gamma(x) = \max_{\xi\in \Xi} \gamma(x, \xi)$, the above reasoning implies that $c(x, \mb P_0)= \bar \gamma(x)$ and that
  \[
    \mb P_0^\infty \left(\hat c_r(x, \hat{\mb P}_T) \ge  c(x, \hat{\mb P}_T)= \bar \gamma(x) \right) = 1 ~ \implies ~ 
    \mb P_0^\infty \left( c(x, \mb P_0) >  \hat c_{r}(x, \hat{\mb P}_T)\right) = 0
    ~ \implies ~ \mb P_0^\infty \left( \hat {\mb P}_T \in {\mc D}(x, \mb P_0) \right) = 0.
  \]
 The probability of being disappointed thus vanishes for all $T\in\mb N$. Hence, it trivially decays at {\em any} rate.

 \paragraph{Case 2:} Assume next that $\mb P_0(\Xi^\star(x))<1$. To prove that the probability of the event $\hat {\mb P}_T\in \mc D(x, \mb P_0)$ decays at a rate of at least $r$, we will first establish the implication
  \begin{equation}
  	\label{eq:disappointment-to-KL}
  	\mb P'\in \bar{\mc D}(x, \mb P_0)\quad \implies \quad \D{\mb P'}{\mb P_0}\ge r.
  \end{equation}

 \paragraph{Case 2a:} Assume that $0<\mb P_0(\Xi^\star(x))<1$. Denote by $\mb U$ the restriction of $\mb P_0$ to $\Xi^\star(x)$, that is, $\mb U(B)=\mb P_0(B\cap\Xi^\star(x))/\mb P_0(\Xi^\star(x))$ for all Borel sets $B\subseteq \Xi$, and define $\mb P(\lambda) = (1-\lambda) \mb P_0+ \lambda \mb U$ for all $\lambda\in[0,1]$. As $\mb P_0(\Xi^\star(x))<1$, one easily verifies that $c(x, \mb P(\lambda))$ is strictly increasing in $\lambda$. The parametric family $\mb P(\lambda)$ now allows us to show that $\D{\mb P'}{\mb P_0}\ge r$ for all $\mb P'\in \bar{\mc D}(x, \mb P_0)$. To this end, assume for the sake of contradiction that~\eqref{eq:disappointment-to-KL} is false and there exists $\mb P'_0\in \bar{\mc D}(x, \mb P_0)$ with $\D{\mb P_0'}{\mb P_0}< r$. Thus, $c(x, \mb P_0) \geq \hat c_r(x, \mb P_0') \geq c(x, \mb P_0)$, where the first inequality holds because $\mb P_0' \in\bar{\mc D}(x, \mb P_0)$, while the second inequality follows from the definition~\eqref{eq:gen:dro} of $\hat c_r$ and the assumption that $\D{\mb P'_0}{\mb P_0}< r$. This reasoning implies that $\mb P_0$ is optimal in~\eqref{eq:gen:dro} for $\mb P'=\mb P'_0$. The assumption $\D{\mb P_0'}{\mb P_0}< r$ further implies that $\mb P_0'$ is absolutely continuous with respect to $\mb P_0$ and consequently also with respect to $\mb P(\lambda)$ for every $\lambda\in [0,1)$. By the definitions of $\mb P(\lambda)$ and $\mb U$, we thus have
  \begin{align*}
    \D{\mb P'_0}{\mb P(\lambda)} & = \int_{\Xi^\star(x)} \log \left( \frac{\d \mb P'_0}{\d (\mb P(\lambda))} \right)\d \mb P'_0 + \int_{\Xi\backslash \Xi^\star(x)} \log \left( \frac{\d \mb P'_0}{\d (\mb P(\lambda))} \right)\d \mb P'_0\\
    &  = \int_{\Xi^\star(x)}  \log \left(\frac{\mb P_0(\Xi^\star(x))}{(1-\lambda)\mb P_0(\Xi^\star(x))+\lambda}\,\frac{\d \mb P'_0}{\d \mb P_0}\right) \d \mb P'_0 + \int_{\Xi\backslash \Xi^\star(x)}  \log \left(\frac{1}{1-\lambda}\,\frac{\d \mb P'_0}{\d \mb P_0}\right) \d \mb P'_0 \\
    & = \D{\mb P'_0}{\mb P_0} + \mb P'_0(\Xi^\star(x))  \log \left(\frac{\mb P_0(\Xi^\star(x))}{(1-\lambda)\mb P_0(\Xi^\star(x))+\lambda} \right) - (1- \mb P'_0(\Xi^\star(x)))\log(1-\lambda),
  \end{align*}
which is continuous in $\lambda\in [0,1)$. The above reasoning implies that there exists $\bar\lambda\in(0,1)$ with $c(x,\mb P(\lambda))>c(x,\mb P_0)$ and $\D{\mb P'_0}{\mb P(\lambda)}\le r$ for all $\lambda\in (0,\bar \lambda]$, which contradicts the optimality of $\mb P_0=\mb P(0)$ in~\eqref{eq:gen:dro} for  $\mb P'=\mb P'_0$. Thus, our assumption was false, and hence \eqref{eq:disappointment-to-KL} follows.

	\paragraph{Case 2b:} Assume now that $\mb P_0(\Xi^\star(x))=0$. In this case we can prove~\eqref{eq:disappointment-to-KL} as in Case~2a. The only differences are that $\mb U$ may now be any distribution on $\Xi^\star(x)$ and that the continuity of $\D{\mb P'_0}{\mb P(\lambda)}$ in $\lambda\in [0,1)$ can now be shown more directly by noting that
  \begin{align*}
    \D{\mb P'_0}{\mb P(\lambda)} & = \int_\Xi \log \left( \frac{\d \mb P'_0}{\d (\mb P(\lambda))} \right)\d \mb P'_0 = \int_\Xi \log \left(\frac{1}{1-\lambda}\,\frac{\d \mb P'_0}{\d \mb P_0}\right) \d \mb P'_0 = \D{\mb P'_0}{\mb P_0} -\log(1-\lambda).
  \end{align*}
All other arguments remain unaffected.

Now that the implication~\eqref{eq:disappointment-to-KL} has been established, we demonstrate that the probability of the event $\hat {\mb P}_T\in \mc D(x, \mb P_0)$ decays at a rate of at least $r$. To this end, we first note that the weak disappointment set $\bar{\mc D}(x, \mb P_0)$ includes the strict disappointment set ${\mc D}(x, \mb P_0)$ and is closed because of the continuity of~$\hat c_r$ established in Proposition~\ref{prop:pred:feasibility:revisited}. The weak LDP upper bound~\eqref{eq:ldp_exponential_rates_ub-continuous} then implies that
  \begin{align*} 
    \limsup_{T\to\infty}\frac 1T \log \mb P_0^\infty \left( \hat {\mb P}_T \in {\mc D}(x, \mb P_0) \right) \leq  -\inf_{\mb P' \in \cl \mc D(x, \mb P_0)} \, \D{\mb P'}{\mb P_0}
    \leq  -\inf_{\mb P' \in \bar{\mc D}(x, \mb P_0)} \, \D{\mb P'}{\mb P_0}\le -r,
  \end{align*}
  where the second inequality holds because $\bar{\mc D}(x, \mb P_0)$ is closed and contains $\mc D(x, \mb P_0)$, while the third inequality follows from~\eqref{eq:disappointment-to-KL}. Thus, the probability of the event $\hat {\mb P}_T\in\mc D(x,\mb P_0)$ decays indeed at a rate of at least~$r$.

As the choice of $x\in X$ and $\mb P_0\in\mc P$ was arbitrary, and as Cases~1 and~2 are exhaustive, $\hat c_r$ is feasible in~\eqref{eq:optimal-predictor}.

In order to show that $\hat c_r$ is strongly optimal in~\eqref{eq:optimal-predictor} when $\epsilon>0$, we can repeat the proof of Theorem~\ref{thm:optimality} almost verbatim with obvious minor modifications (most notably, there is no need to construct $\mb P_2$).
   \hfill $\square$
  \endproof

  \proof{Proof of Theorem \ref{thm:presc:feasibility-revisited}.}
We first establish feasibility for $r\ge 0$ and $\epsilon>0$. From Proposition~\ref{prop:pred:feasibility:revisited} and the subsequent discussion we know that $\hat c_r$ is continuous and that $\hat x_r$ can be chosen to be quasi-continuous, which implies that $(\hat c_r+\epsilon,\hat x_r)\in\mc X$. To show that the out-of-sample disappointment of $(\hat c_r+\epsilon, \hat x_r)$ decays sufficiently fast, we fix $\mb P_0\in\mc P$ and define
 \[
    \mc D_\epsilon(\mb P_0)= \left\{ \mb P'\in\mc P : c(\hat x_r(\mb P'), \mb P_0) >  \hat c_{r}(\hat x_r(\mb P'), \mb P')  + \epsilon\right\}
\]
as the corresponding disappointment sets. For any $x\in X$ and $\mb P_0\in\mc P$ we further~define
  \[
    \mc D_\epsilon(x, \mb P_0)= \left\{ \mb P'\in\mc P : c(x, \mb P_0) >  \hat c_{r}(x, \mb P') + \epsilon \right\}\quad \text{and}\quad
    \bar{\mc D}_\epsilon(x, \mb P) = \left\{ \mb P' \in\mc P: c(x, \mb P) \ge  \hat c_{r}(x, \mb P') + \epsilon \right\}
  \]
as the corresponding strict and weak {\em decision-dependent} disappointment sets. In order to show that the probability of the event $\hat {\mb P}_T\in \mc D_\epsilon (x, \mb P_0)$ decays at a rate of at least $r$, we observe that
\begin{align*}
  \mc D_\epsilon (\mb P_0) \subseteq \bigcup_{x\in X} \mc D_\epsilon(x, \mb P_0) \subseteq \bigcup_{x\in X} \bar{\mc D}_\epsilon(x, \mb P_0)= \bar{\mc D}_\epsilon(\mb P_0),
\end{align*}
where 
\[
	\bar{\mc D}_\epsilon(\mb P_0)=\set{\mb P'\in\mc P}{\max_{x\in X} c(x, \mb P_0)-\hat c_r(x, \mb P')\geq \epsilon}.
\]
The proof of Theorem~\ref{thm:pred:feasibility-revisited} immediately implies that $\bar {\mc D}_\epsilon(x, \mb P_0)$ is closed. We will now argue that $\bar{\mc D}_\epsilon(\mb P_0)$ is also closed. As the model-based predictor $c$ and the distributionally robust predictor $\hat c_r$ are both jointly continuous in $x$ and $\mb P'$ and as the feasible set $X$ is compact, the maximum theorem by \citet[pp.~115--116]{berge1963topological} implies that the function $\max_{x\in X} c(x, \mb P_0)-c_r(x, \mb P')$ is continuous in $\mb P'$ for any fixed $\mb P_0$. Thus, $\bar {\mc D}_\epsilon(\mb P_0)$ is closed as a superlevel set of a continuous function. As $\mc D_\epsilon (\mb P_0) \subseteq \bar{\mc D}_\epsilon(\mb P_0)$, this implies that $\cl \mc D_\epsilon (\mb P_0) \subseteq \bar{\mc D}_\epsilon(\mb P_0)$.

Next, we will establish the implication
  \begin{equation}
  	\label{eq:disappointment-to-KL-epsilon}
  	\mb P'\in \bar{\mc D}_\epsilon(x, \mb P_0)\quad \implies \quad \D{\mb P'}{\mb P_0}\ge r
  \end{equation} 
for any $\epsilon>0$. Assume first that $\mb P_0(\Xi^\star(x))=1$, where $\Xi^\star(x)= \arg \max_{\xi\in \Xi} \gamma(x, \xi)$ stands as usual for the (compact) set of worst-case scenarios. Then, for any $\mb P'\in\mc P$ with $\D{\mb P'}{\mb P_0}< r$ we have
  \[
	\mb P'\ll \mb P_0\quad \implies\quad \mb P'(\Xi^\star(x))=1\quad 
	\implies \quad c(x,\mb P_0)=\hat c_r(x,\mb P') \quad \implies\quad \mb P'\notin \bar {\mc D}_\epsilon(x,\mb P_0),
\]
where the first implication holds because the support of $\mb P_0$ is a assumed to be a subset of $\Xi^\star(x)$ and because the support of any distribution $\mb P'$ that is absolutely continuous with respect to $\mb P_0$ must be contained in the support of $\mb P_0$. The second implication follows from the observation that both $c(x,\mb P_0)$ and $\hat c_r(x,\mb P')$ must evaluate to the worst-case cost $\bar \gamma(x)=\max_{\xi\in \Xi} \gamma(x,\xi)$ because both $\mb P_0$ and $\mb P'$ are supported on the set of worst-case scenarios $\Xi^\star(x)$. Thus, $\D{\mb P'}{\mb P_0}< r$ implies $\mb P'\notin \bar {\mc D}_\epsilon(x,\mb P_0)$, whereby \eqref{eq:disappointment-to-KL-epsilon} follows by contraposition.

Assume next that $\mb P_0(\Xi^\star(x))<1$. Then~\eqref{eq:disappointment-to-KL-epsilon} is an immediate consequence of the stronger implication~\eqref{eq:disappointment-to-KL} derived in the proof of Theorem~\ref{thm:pred:feasibility-revisited}.

In summary, we thus find
\begin{equation*}
\begin{aligned}
  \limsup_{T\to\infty}\frac 1T \log \mb P_0^\infty \left( \hat {\mb P}_T \in \mc D_\epsilon(\mb P_0) \right) \leq -\inf_{\mb P'\in \bar{\mc D}_\epsilon(\mb P_0)}\D{\mb P'}{\mb P_0} = -\inf_{x\in X} \inf_{\mb P'\in \bar{\mc D}_\epsilon(x, \mb P_0)} \D{\mb P'}{\mb P_0} \leq -r,
\end{aligned}
\end{equation*}
where the first inequality follows from the weak LDP upper bound \eqref{eq:ldp_exponential_rates_ub-continuous} and the inclusion $\cl \mc D_\epsilon (\mb P_0) \subseteq \bar{\mc D}_\epsilon(\mb P_0)$. The equality exploits the definition of $\bar{\mc D}_\epsilon(\mb P_0)$, and the second inequality follows from the inclusion~\eqref{eq:disappointment-to-KL-epsilon}, which holds for any $\epsilon>0$. As $\mb P_0\in\mc P$ and $\epsilon>0$ were chosen arbitrarily, $(\hat c_r+\epsilon,\hat x_r)$ is thus feasible in~\eqref{eq:optimal-prescriptor}. 

To show that $(\hat c_r,\hat x_r)$ is preferred to any feasible solution in~\eqref{eq:optimal-prescriptor} when $\epsilon>0$, we can repeat the proof of Theorem~\ref{thm:optimality_prescriptor} almost verbatim with obvious minor modifications ({\em e.g.}, there is no need to construct $\mb P_2$).
 \hfill $\square$
\endproof

\end{APPENDICES}

\newpage

\bibliographystyle{plainnat}
\bibliography{references}
 
%%%%%%%%%%%%%%%%%
\end{document}